%% PartI.TeX, plain TeX file.
%% Preliminary version of July 30, 2004
%% not to be published or distributed.

\magnification=\magstep1

%-----Paper Format---------
\hsize=14cm
\vsize=18cm
\hoffset=-0.3cm
\voffset=0cm
\footline={\hss{\vbox to 1cm{\vfil\hbox{\rm\folio}}}\hss}

%-----AMS Symbols-------
\input amssym.def
\input amssym.tex

%-----Special Fonts-----

%\font\title=cmr12 at 16pt
\font\rr=cmcsc10 scaled 1200
\font\teneusm=eusm10
\font\seveneusm=eusm7
\font\fiveeusm=eusm5
\newfam\eusmfam
\def\eusm{\fam\eusmfam\teneusm}
\textfont\eusmfam=\teneusm
\scriptfont\eusmfam=\seveneusm
\scriptscriptfont\eusmfam=\fiveeusm

\font\teneufm=eufm10
\font\seveneufm=eufm7
\font\fiveeufm=eufm5
\newfam\eufmfam

\textfont\eufmfam=\teneufm
\scriptfont\eufmfam=\seveneufm
\scriptscriptfont\eufmfam=\fiveeufm

%----Small Macros----
\def\varGamma{{\mit \Gamma}}
\def\Re{{\rm Re}\,}
\def\Im{{\rm Im}\,}
\def\txt#1{{\textstyle{#1}}}
\def\scr#1{{\scriptstyle{#1}}}
\def\r#1{{\rm #1}}
\def\B#1{{\Bbb #1}}
\def\e#1{{\eusm #1}}
\def\b#1{{\bf #1}}

\footnote{{}}{\hskip -2 mm 2001 {\it Mathematics Subject
Classification}.
11M06, 11F72, 11F66

}
\footline={\hss{\vbox to 2cm{\vfil\hbox{\rm\folio}}}\hss}
\nopagenumbers
   
   \font\bb=cmcsc10

\def\txt#1{{\textstyle{#1}}}
\baselineskip=13pt
\def\hf{{\textstyle{1\over2}}}
\def\a{\alpha}
\def\d{{\,\rm d}}
\def\eps{\varepsilon}

\def\G{\Gamma}
\def\k{\kappa}
\def\s{\sigma}

\def\={\;=\;}
\def\zx{\zeta(\hf+ix)}
\def\zt{\zeta(\hf+it)}
\def\zr{\zeta(\hf+ir)}

\def\no{\noindent}
  \def\s{\sigma}
\def\z{\zeta}
\def\b{\beta}

\def\no{\noindent} \def\d{\,{\rm d}} 
  \def\s{\sigma}
      \def\G{\Gamma}

\font\teneufm=eufm10
\font\seveneufm=eufm7
\font\fiveeufm=eufm5
\newfam\eufmfam
\textfont\eufmfam=\teneufm
\scriptfont\eufmfam=\seveneufm
\scriptscriptfont\eufmfam=\fiveeufm
\def\mathfrak#1{{\fam\eufmfam\relax#1}}

\font\tenmsb=msbm10

\newfam\msbfam
       \textfont\msbfam=\tenmsb
%      \scriptfont\msbfam=\sevenmsb
%      \scriptscriptfont\msbfam=\fivemsb
\def\Bbb#1{{\fam\msbfam #1}}

\def \NN {\Bbb N}
\def \CC {\Bbb C}

\def \RR {\Bbb R}
\def \ZZ {\Bbb Z}

   \def\rightheadline{{\hfil{\sevenrm
   The fourth moment off the critical line  }\hfil\sevenrm\folio}}

   \def\leftheadline{{\sevenrm\folio\hfil{\sevenrm
   A. Ivi\'c and Y. Motohashi }\hfil}}
   \def\emptyheadline{\hfil}
   \headline{\ifnum\pageno=1 \emptyheadline\else
   \ifodd\pageno \rightheadline \else \leftheadline\fi\fi}

\def \vG {{\mit \Gamma}}
%-------Text---------
\centerline{\rr The Moments of the Riemann Zeta-Function.}
\bigskip
\centerline{\rr Part I: The fourth moment off the critical line}
\medskip
\centerline{{\bb Aleksandar Ivi\'c and Yoichi Motohashi}}
\bigskip
\centerline{\bb Abstract}

\bigskip\no
In this paper, the first part of a larger work, we prove the
spectral decomposition of
$$
\int_{-\infty}^\infty|\zeta(\s+it)|^4g(t)\d t\qquad(\hf < \s < 1\;\; {\rm {fixed}}),
$$
where $g(t)$ is a suitable weight function of fast decay. This is used to
obtain estimates and omega results for the function
$$\eqalign{
E_2(T,\s) &: =\int_0^T|\z(\s+it)|^4\d t -
{\z^4(2\s)\over\z(4\s)}T -{T\over3-4\s}{\left({T\over2\pi}
\right)}^{2-4\s}{\z^4(2-2\s)\over\z(4-4\s)}\cr&
\,\,- T^{2-2\s}(a_0(\s) + a_1(\s)\log T + a_2(\s)\log^2T),\cr}
$$
the error term in the asymptotic formula for the fourth moment
of $|\z(\s+it)|$.

\vskip 1cm
\centerline{\bb 1. Introduction}
\bigskip

Power moments of the Riemann zeta-function $\z(s)$ are one of the
central objects in the theory of $\z(s)$, with many important
applications.
Although the main interest is in the moments on the ``critical line"
$\Re s = \hf$, the moments when $s$ lies in the
``critical strip" $\hf < \Re s < 1$, or
``off" the critical line, are also of great interest. There exist
extensive results on the second and fourth moments on the critical
line,
the only ones that so far can be treated unconditionally, and where
asymptotic
formulas have been obtained.
A comprehensive review on mean square results for $\z(s)$
is given by Matsumoto [Ma], where further references may be found. Some
of the relevant works on the fourth moment of
$|\zt|$ are [I1], [I2], [I4]--[I8], [IM1]--[IM3], [IJM], [Mo2]--[Mo6],
where also the interested reader may find
further references. The
aim of this paper is to treat the fourth moment off the critical line.
The only
works that seem to have explicitly  dealt with
this
subject are [I6], [K1] and [K2]. Thus it appears that the time is ripe
for an
extensive account of this subject, which we hope that the present work
will
provide.

\medskip
  The main object of our study is the weighted integral
$$
\e{L}(g;\s,\tau) \;:=\; \int_{-\infty}^\infty
\left|\zeta\left(\s+it\right)\right|^2|\zeta(\tau+it)|^2g(t)\d t,
\leqno(1.1)
$$
where  $\s,\tau$ are given constants satisfying
$$
\hf \le \s \le \tau\quad(\s \ne 1,\; \tau\ne 1).\leqno(1.2)
$$
\medskip
\noindent
{\it The basic assumption on the weight $g$ is\/}:
The even function $g(t)$ takes real values on the real axis, and there
exists a large positive constant $A$ such that $g(t)$ is regular and
$g(t) = O((|t|+1)^{-A})$ in the horizontal strip $|\Im{t}|\le{A}$.
\medskip
We shall obtain the spectral decomposition of $\e{L}(g;\s,\tau)$
(see Section 3) by the method used by the second author in the
case of $\e{L}(g;\hf,\hf)$ (see [Mo2], [Mo6]). This decomposition,
which
is in fact an exact identity, will contain, among other things, the
function
$$
g^*(\xi)=\int_{-\infty}^\infty g(t){\rm e}^{-i\xi{t}}\d
t\qquad(\xi\in\RR),\leqno(1.3)
$$
namely the Fourier transform of $g$. Note that, since $g$ is even,
$$
g^*(\xi) = \int_{-\infty}^\infty g(t){\rm e}^{i\xi{t}}\d t
= \int_{-\infty}^\infty g(t)\cos(\xi{t})\d t = g_c(\xi),\leqno(1.4)
$$
where $g_c$ is the {\it cosine Fourier transform} of $g$.
The function $\e{L}(g;\s,\tau)$ is, with an appropriate choice of
the weight  $g$, {\it the local object} which
after the integration over a suitable parameter contained in $g$ will
lead to
the asymptotic evaluation of {\it the global object}
$$
\int_0^T|\z(\s+it)|^2|\z(\tau+it)|^2\d t,\leqno(1.5)
$$
provided that (1.2) holds. A good choice of $g$ will
entail rapid decay of $g_c$, which will facilitate handling
of the quantities that will appear in the spectral decomposition.

\medskip
It is clear that (1.5) is not
interesting when $\s > 1, \tau > 1$, in which case the
zeta-values in question are represented by absolutely convergent
series which may be readily integrated termwise.
The special cases of interest of (1.1) and (1.5) are

\medskip
a) $\s = \tau = \hf$. This is the classical case of the fourth
moment of $\z(s)$ on the critical line, and probably
the most important case. It was obtained by the second author
[Mo2], and is extensively discussed in  [I2] and [Mo6].

\medskip
b) $\hf < \s = \tau < 1$. This is the  case of the fourth
moment of $\z(s)$ off the critical line. As already mentioned,
this is discussed by the first author in [I6]
and by A. Ka\v c\.enas [K1], [K2]. The formula for the fourth
moment reads (when $\hf < \s < 1$ is fixed)
$$
\leqalignno{\int_0^T|\z(\s+it)|^4\d t &=
{\z^4(2\s)\over\z(4\s)}T +
{T\over3-4\s}{\left({T\over2\pi}\right)}^{2-4\s}{\z^4(2-2\s)\over\z(4
-4\s)}
&(1.6)\cr&
+ T^{2-2\s}(a_0(\s) + a_1(\s)\log T + a_2(\s)\log^2T) +
E_2(T,\s),\cr}
$$
where $E_2(T,\s)$ is the error term, and the $a_j(\s)$'s are constants
which may be explicitly evaluated. When $\s\to\hf+0$ the
function $E_2(T,\s)$ tends to $E_2(T,\hf)\equiv E_2(T)$, the
error term in the asymptotic formula for $\int_0^T|\zt|^4\d t$.

\medskip
In [I6] only a sketch of the spectral
decomposition of the fourth moment off the critical line,
due to the second author, was given. Here we are
going to give a rigorous proof of the spectral decomposition in
question and to recover and extend the results given in [I6].
The works of Ka\v c\.enas contain an explicit
evaluation of the
main term in the asymptotic formulas for the fourth moment off the
critical
line, but the estimates for the error term are weaker than those given
in
[I6]. We also remark that the cases a) and b) have their analogues
(mean
squares) for automorphic $L$--functions (see [Mo3], [Mo7]). The
fourth moment of $\z(s)$ off the critical line has its analogue in
the mean square off the critical line. In this case, which is less
difficult to deal with than the present case, the formula reads
$$
\int_0^T|\z(\s+it)|^2\d t = \z(2\s)T +
(2\pi)^{2\s-1}\,{\z(2-2\s)\over2-2\s}
T^{2-2\s} + E_1(T,\s)\quad(\hf < \s < 1),\leqno(1.7)
$$
where $E_1(T,\s)$ in (1.7) represents the error term, and the
notational
analogy between $E_1(T,\s)$ and $E_2(T,\s)$ is obvious. As we
already mentioned, [Ma] represents a comprehensive survey of
results on $E_1(T,\s)$.

\medskip
c) $\s = \hf, \hf < \tau < 1$. This case, which does not seeem to
have been treated in the literature before, may be thought of as
a ``hybrid mean value".

\medskip
d) $\s = \hf, \tau > 1$. This case is an extension of c).
When $\tau$ is large, it is of interest because then $\e{L}(g;\hf,\tau)$
becomes a multiple of
$$
\int_{-\infty}^\infty
\left|\zeta\left(\hf+it\right)\right|^2g(t)\d t,
$$
and provides the mean
(at least theoretically)
to estimate $\zt$ pointwise, which is a fundamental
problem in the theory of $\z(s)$.

\medskip
In this paper, which is Part I of the whole work, we shall treat the
case
b) above. To avoid excessive length, the cases c) and d) will be
treated
in Part II. The plan of the present paper is as follows. The
formulation of the spectral decomposition of
${\e{L}}(g;\tau,\tau)$, when $\hf < \tau < 1,$ will be given
in Section 2. Although the proof has many analogies with the proof
of the second author for the case of ${\e{L}}(g;\hf,\hf)$, there
are also many detours, and the complete, rigorous proof is given
in Section 3. In the result  $\a_j H^2_j(\hf)H_j(\tau)$
and $\a_j H_j(\hf)H_j^2(\tau)$ (in Part II)
appear, and the asymptotic
evaluation of sums of these quantities over $\k_j\le K$
is carried out in Section
4. The detailed asymptotic evaluation of the function $\Lambda$,
appearing in the spectral decomposition of ${\e{L}}(g;\tau,\tau)$
with the Gaussian weight function, is contained in Section 5.
The explicit formulas for  ${\e{L}}(g;\tau,\tau)$ and its integral
are presented in Section 6. They
are necessary in order to obtain results on the error term
$E_2(T,\s)$, which is done in Section 7 and Section 8. The
notation used throughout the paper is, whenever possible,
standard. We have used the letter $\tau$ occasionally where one
would commonly used $\s$ (as in the notation for $E_2(T,\s)$).
This was done to avoid possible confusion with the real part of
the complex variable $s$, especially in Section 3.

\bigskip
\centerline{\bb 2. Spectral decomposition of the fourth moment
-- notation and results}
\bigskip
In this section we introduce the necessary notation for the spectral
decomposition of ${\e{L}}(g;\tau,\tau)$, the weighted fourth moment
off the critical line. We also present Theorem 1, which will give the
desired decomposition, but postpone the proof for Section 3. The
notation
used throughout is standard, to be found e.g., in the second author's
monograph [Mo6], and for this reason we shall be relatively brief.

Let $\,\{\lambda_j = \kappa_j^2 + {1\over4}\} \,\cup\, \{0\}\,$ be the
discrete spectrum of the hyperbolic Laplacian
$$
\Delta=-y^2\left({\left({\partial\over\partial x}\right)}^2 +
{\left({\partial\over\partial y}\right)}^2\right)
$$
acting over the Hilbert space composed of all
$\vG$-automorphic functions which are square integrable with
respect to the hyperbolic measure, where
$$
\vG \;\cong\; SL(2,\,\ZZ)/\{+1,-1\}.
$$
  Let $\{\psi_j\}$ be a maximal
orthonormal system in this space such that
$\Delta\psi_j=\lambda_j\psi_j$ for each $j\ge1$ and
$T(n)\psi_j=t_j(n)\psi_j$  for each integer $n\in\NN$, where
$$
\bigl(T(n)f\bigr)(z)\;=\;{1\over\sqrt{n}}\sum_{ad=n}
\,\sum_{b=1}^df\left({az+b\over d}\right)
$$
is the Hecke operator. We shall further assume that
$\psi_j(-\bar{z})=\eps_j\psi_j(z)$ with the parity sign $\eps_j=\pm1$. We
then define ($s = \s + it$ will denote a complex variable)
$$
H_j(s)\;=\;\sum_{n=1}^\infty t_j(n)n^{-s}\qquad(\s > 1),\leqno(2.1)
$$
which denotes the Hecke series associated with
$\psi_j(z)$, and which can be continued to an entire function.
As usual  we put
$$
\a_j = |\rho_j(1)|^2(\cosh\pi\kappa_j)^{-1},\leqno(2.2)
$$
where $\rho_j(1)$ is the first Fourier coefficient of
$\psi_j(z)$. The holomorphic counterparts $\a_{j,k}$ and
$$
H_{j,k}(s) := \sum_{n=1}^\infty t_{j,k}(n)n^{-s}\qquad(\Re s > 1)
$$
of (2.2) and (2.1), respectively, are defined in [Mo6, Chapter
3]; as to
$\vartheta(2k)$ in $(2.7)$ below see Section 2.2 there. Now we can
formulate

\medskip
THEOREM 1.
{\it Let ${1\over2}<\tau<1$ be fixed, and let $g$ satisfy the
basic assumption. Then we have} ({\it cf.} (1.1))
$$
{\e{L}}(g;\tau,\tau)
= \big\{{\e{Z}}_r+{\e{Z}}_d+{\e{Z}}_c +{\e{Z}}_h\big\}(\tau,g),
\leqno(2.3)
$$
{\it where
$$
\leqalignno{
&{\e{Z}}_r(\tau,g)=M(\r{p}_\tau;g)&(2.4)\cr
-&8\pi\zeta(2\tau-1)^2\Re\left\{\left(c_E-{\zeta'\over\zeta}(2\tau-1)
\right)g((\tau-1)i)+{1\over2}ig'((\tau-1)i)\right\},\cr
}
$$
with the function $M$ being defined by $(3.65)$, $(3.88)$, and $(3.92)$
according
as ${1\over2}<\tau<{3\over4}$, $\tau={3\over4}$, and
${3\over4}<\tau<1$,
respectively. Further we have
$$
{\e{Z}}_d(\tau,g)=\sum_{j=1}^\infty
\alpha_jH_j^2(\hf)H_j\left(2\tau-\hf\right)
\Lambda(\kappa_j;\tau,g),\leqno(2.5)
$$
$$
{\e{Z}}_c(\tau,g)={1\over\pi}\int_{-\infty}^\infty
{{|\zeta({1\over2}+ir)|^4|\zeta(2\tau-{1\over2}+ir)|^2}
\over{|\zeta(1+2ir)|^2}}\Lambda(r;\tau,g)\d r,\leqno(2.6)
$$
$$
{\e{Z}}_h(\tau,g)=\sum_{k=1}^\infty\sum_{j=1}^{\vartheta(2k)}
\alpha_{j,2k}H_{j,2k}^2
(\hf)H_{j,2k}(2\tau-\hf)
\Lambda\left((\hf-2k)i\,;\tau,g\right).\leqno(2.7)
$$
Here $c_E = -\G'(1)$ is Euler's constant, and
$$
\leqalignno{
&\Lambda(r;\tau,g)=\int_0^\infty
(y(1+y))^{-\tau}g_c\left(\log\left(1+{1\over y}\right)\right)&(2.8)\cr
&\times{\Re}\left\{y^{-{1\over2}-
ir}\left(1+{i\over{\sinh(\pi{r})}}\right)
{{\Gamma({1\over2}+ir)^2}\over{
\Gamma(1+2ir)}} F\left(\hf+ir,\hf+ir;1+2ir;-{1\over y}\right)
\right\}\d y\cr
}
$$
with the hypergeometric function $F$}.

\medskip

The above spectral decomposition is analogous to the spectral
decomposition of the function ${\e{L}}(g;\hf,\hf)$, given as [Mo6,
Theorem 4.2].
It is in fact an exact identity, relating the original object
(weighted integral of the fourth moment) to various objects from
spectral theory, hence the terminology ``spectral decomposition".
The notation is also analogous to the one used
in [Mo6, Theorem 4.2], as much as possible. The notation
$M(\r{p}_\tau;g)$ refers to the ``main term", since suitable
integration
of this term will lead to the main term for the fourth moment of
$|\z(\s+it)|$ itself (see (1.6)). Likewise, the notation
${\e{Z}}_r, {\e{Z}}_d, {\e{Z}}_c, {\e{Z}}_h$ refers to
``residual", ``discrete", ``continuous" and ``holomorphic" parts,
respectively. As we just mentioned, the term $M(\r{p}_\tau;g)$,
contained in ${\e{Z}}_r$, will eventually contribute to the main
term, while the remaining terms will contribute to the error terms.
Of these, the most difficult (major) contribution, like in the case
of ${\e{L}}(g;\hf,\hf)$, will come from
${\e{Z}}_d$.

\smallskip
An important feature of the above formula is the appearance of the
oscillatory integral $\Lambda(r;\tau,g)$ which containins the hypergeometic
function. We recall here that, for $|z| < 1$, one  defines the
hypergeometric function
$$
\leqalignno{
F(\a,\b;\gamma;z) &= \sum_{k=0}^\infty
{(\a)_k(\b)_k\over(\gamma)_kk!}z^k&(2.9)\cr& \= 1 +
\sum_{k=1}^\infty{\a(\a+1)\ldots(\a+k-1)\b(\b+1)\ldots(\b+k-1)\over
\gamma(\gamma+1)\ldots(\gamma+k-1)k!}z^k.\cr}
$$
Analytic continuation and other properties of $F(\a,\b;\gamma;z)$
are treated e.g., by N.N. Lebedev [L].

\bigskip
\centerline{\bb 3. Proof of the spectral decomposition
for the fourth moment}
\bigskip
This section contains the proof of Theorem 1; we assume throughout that
$\hf < \tau < 1$ is fixed, and that the {\it basic assumption}
on $g$ holds. Our argument is a reworking of [Mo6, Chapter 4];
thus we could mention specific changes only. However, that would make
the later
part of our discussion hard to comprehend, since as has been mentioned
above there
are many sensitive detours peculiar to our new situation that begins in
fact
at $(3.7)$ below.

Let first
$$
\leqalignno{
\tilde{g}(s,\lambda)&=\int_0^\infty
y^{s-1}(1+y)^{-\lambda}g^*(\log(1+y))\d y
&(3.1)\cr
&=\Gamma(s)\int_{-\infty+Ai}^{\infty+Ai}{\Gamma(\lambda-it-s)\over
\Gamma(\lambda-it)}g(t)\d t,
}
$$
where $g^*$ is defined by (1.3). We begin with the
analogue of [Mo6, Lemma 4.1], namely
\medskip
{\bf Lemma 1.}
{\it The function $\tilde{g}(s,\lambda)/\Gamma(s)$ continues
holomorphically to the domain
$$
|\Re{s}|\le {1\over3}A,\quad |\Re\lambda|\le
{1\over3}A\,;\leqno(3.2)
$$
and there we have
$$
\tilde{g}(s,\lambda)\ll |s|^{-{1\over2}A},
\leqno(3.3)
$$
when $s$ tends to infinity while $\lambda$ remains bounded.
}

\medskip
Let now ${\cal D}_+$ and ${\cal D}_-$ be the domains of ${\Bbb C}^4$
where all
four variables have real parts larger than and less than one,
respectively. We set, for $(u,v,w,z)\in {\cal D}_+$,
$$
\e{J}(u,v,w,z;g) :=\int_{-\infty}^\infty\zeta(u+it)
\zeta(v+it)\zeta(w-it)\zeta(z-it)g(t)\d t.
\leqno(3.4)
$$
Moving the path upwards appropriately, we see that ${\e J}$ is a
meromorphic function over the domain
$$
{\cal B}\;=\;\{\,(u,v,w,z)\in{\Bbb
C}^4\,:\,|u|,\,|v|,\,|w|,\,|z|<B\,\},
\leqno(3.5)
$$
where $B=cA$ with $0<c<1$ is supposed to be sufficiently large. Then,
taking $(u,v,w,z)$ in ${\cal D}_-\cap{\cal B}$, we get the following
meromorphic continuation of ${\e J}$ to ${\cal D}_-\cap{\cal B}$:
$$
\leqalignno{
{\e J}(u,v,&w,z;g)=\int_{-\infty}^\infty\zeta(u+it)
\zeta(v+it)\zeta(w-it)\zeta(z-it)g(t)\d t&(3.6)\cr
&+2\pi\zeta(v-u+1)\zeta(u+w-1)\zeta(u+z-1)g((u-1)i)\cr
&+2\pi\zeta(u-v+1)\zeta(v+w-1)\zeta(v+z-1)g((v-1)i)\cr
&+2\pi\zeta(z-w+1)\zeta(u+w-1)\zeta(v+w-1)g((1-w)i)\cr
&+2\pi\zeta(w-z+1)\zeta(u+z-1)\zeta(v+z-1)g((1-z)i).\cr
}
$$
\medskip
{\bf Lemma 2.}
{\it The function ${\e J}$ is regular at the point
$\r{p}_\tau \;:=\;(\tau,\tau,\tau,\tau)$,
and we have
$$
\leqalignno{
&{\e{L}}(g;\tau,\tau)={\e{J}}(\r{p}_\tau;g)&(3.7)\cr
&-8\pi\zeta(2\tau-1)^2\Re\left\{\left(c_E-{\zeta'\over\zeta}(2\tau-1)
\right)g((\tau-1)i)+{1\over2}ig'((\tau-1)i)\right\},
}
$$
where $c_E = -\G'(1)$ is  Euler's constant.
}
\medskip
{\bf Proof.\/} On the right side of $(3.6)$ the integral is obviously
regular
throughout
${\cal D}_-\cap{\cal B}$. To see the regularity at $\r{p}_\tau$ of
the sum of other terms, we need only to replace the factors
$\zeta(v-u+1)$,
$\zeta(u-v+1)$, $\zeta(z-w+1)$ and $\zeta(w-z+1)$ by their Laurent
expansions. For example, its value at $\r{p}_\tau$ is
$$
\leqalignno{
&4\pi\zeta(2\tau-1)^2c_E\{g((\tau-1)i)+g((1-\tau)i)\}&(3.8)\cr
-&2\pi\Big[{\partial\over{\partial{v}}}\{\zeta(v+w-1)\zeta(v+z-1)g((v
-1)i)\}
\Big]_{\r{p}_\tau}\cr
-&2\pi\Big[{\partial\over{\partial{z}}}\{\zeta(u+z-1)\zeta(v+z-1)g((1-
z)i)\}
\Big]_{\r{p}_\tau}\cr
=\,&4\pi\{\zeta(2\tau-1)^2c_E-\zeta(2\tau-1)\zeta'(2\tau-1)\}
\{g((\tau-1)i)+g((1-\tau)i)\}\cr
+\,&2\pi\zeta(2\tau-1)^2i\{g'((\tau-1)i)-g'((1-\tau)i)\}.\cr
}
$$
Next, in ${\cal D}_+$ we have
$$
\leqalignno{
{\e J}(u,v,w,z;g)&=\sum_{k,l,m,n=1}^\infty
k^{-u}l^{-v}m^{-w}n^{-z}g^*(\log(mn)/(kl))&(3.9)\cr
&={\e J}_0(u,v,w,z;g)
+{\e J}_1(u,v,w,z;g)+\overline{{\e J}_1(\overline{w},\overline{z},
\overline{u},\overline{v};g)},
}
$$
where ${\e J}_0$ and ${\e J}_1$ correspond to the parts with $kl=mn$
and
$kl<mn$, respectively. We have
$$
{\e J}_0(u,v,w,z;g)= g^*(0)\zeta(u+w)\zeta(u+z)\zeta(v+w)
\zeta(v+z)/\zeta(u+v+w+z),\leqno(3.10)
$$
and
$$
{\e J}_1(u,v,w,z;g)={1\over{2\pi{i}}}\sum_{m,n=1}^\infty
{{\sigma_{u-v}(m)\sigma_{w-z}(m+n)}\over{m^{u+w}}}
\int_{(2)}\tilde{g}(s,w)(m/n)^s\d s,
\leqno(3.11)
$$
where $\int_{(c)}\cdots \d s$ denotes integration over the line
$\Re s = c$, and $\s_a(n) = \sum_{d|n}d^a$.
One may deal with this double sum in two ways: either by using
the Ramanujan expansion of the function $\sigma_{w-z}(m+n)$,
or by embedding ${\e J}_1$  in values of a Poincar\'e
series on $\varGamma\backslash \r{PSL}_2(\B{R})$. Here the
first method is employed, and we shall follow [Mo6, Chapter 4].
As to the second method, see [BM]. It should be remarked that the
latter
dispenses with the spectral theory of sums of Kloosterman sums that
plays a predominant r\^ole in the former. Also it should be added in
this context
that Theorem 1 above could be formulated solely in terms of the
$\varGamma$-automorphic representations of
$\r{PSL}_2(\B{R})$.
\medskip
{\bf Lemma 3.}
{\it The function ${\e J}_1(u,v,w,z;g)$ can be continued
meromorphically to the domain
$$
{\cal E}:=\{(u,v,w,z)\in{\cal B}: \Re(u+w)< {\txt{1\over3}}B,
\Re(v+w)< {\txt{1\over3}}B, \Re(u+v+w+z)>3B\},\leqno(3.12)
$$
and in $\cal E$ we have the decomposition
$$
{\e J}_1(u,v,w,z)={\e J}_2(u,v,w,z)+{\e J}_3^{+}(u,v,w,z)
+{\e J}_3^{-}(u,v,w,z).
\leqno(3.13)
$$
Here
$$
\leqalignno{
{\e J}_2&(u,v,w,z;g)
:=\tilde{g}(u+w-1)\zeta(v+z)\zeta(u+w-1)&(3.14)\cr
&\qquad\times\zeta(z-w+1)\zeta(v-u+1)/\zeta(v+z-u-w+2),\cr
&+\tilde{g}(v+w-1)\zeta(u+z)\zeta(v+w-1)\cr
&\qquad\times\zeta(z-w+1)\zeta(u-v+1)/\zeta(u+z-v-w+2)
}
$$
and
$$
\leqalignno{
&{\e J}_3^\pm(u,v,w,z;g):=2(2\pi)^{w-z-1}\zeta(z-w+1)&(3.15)
\cr
&\times\sum_{m,n=1}^\infty
m^{{1\over2}(1-u-v-w-z)}n^{{1\over2}(u+w-v-z-1)}
\sigma_{v-u}(n)K_\pm(m,n;u,v,w,z;g),
}
$$
where
$$
K_\pm(m,n;u,v,w,z;g)=\sum_{l=1}^\infty{1\over{l}}S(m,\pm{n};l)
\varphi_\pm\Big({{4\pi}\over{l}}\sqrt{mn};u,v,w,z;g\Big)
$$
with $S(a,b;c) = \sum_{1\le n\le c,(n,c)=1,n\bar {n}\equiv1({\rm
mod}\,c)}
\exp\left(2\pi i\bigl({an+b\bar {n}\over c}\bigr)\right)$ a Kloosterman
sum, and
$$
\leqalignno{
\varphi_+&(x;u,v,w,z;g):={1\over{2\pi{i}}}\cos(\hf(u-v)\pi)
\int_{(B)}\left({x\over2}\right)^{u+v+w+z-1-2s}&(3.16)\cr
&\times\Gamma(s+1-u-w)\Gamma(s+1-v-w)\tilde{g}(s,w)\d s,
}
$$
$$
\leqalignno{
\varphi_-(x;u,v,w,z;g)&:=-{1\over{2\pi{i}}}\int\limits_{(B)}
\left({x\over2}\right)^{u+v+w+z-1-2s}\cos(\pi(w+
\hf(u+v)-s))&(3.17)\cr
&\times\Gamma(s+1-u-w)\Gamma(s+1-v-w)\tilde{g}(s,w)\d s.
}
$$
}
\medskip
The Kloosterman--Spectral sum formula of N.V. Kuznetsov (see [Mo6])
yields, with the
standard notation from the spectral theory of the Fourier coefficients
of
modular cusp forms, that
$$
\leqalignno{
&K_+(m,n;\,u,v,w,z;g):=\sum_{j=1}^\infty\alpha_jt_j(m)t_j(n)
(\varphi_+)^+(\kappa_j;u,v,w,z;g)&(3.18)\cr
&+{1\over\pi}\int_{-\infty}^\infty{{\sigma_{2ir}(m)
\sigma_{2ir}(n)}\over{(mn)^{ir}|\zeta(1+2ir)|^2}}
(\varphi_+)^+(r;u,v,w,z;g)\d r\cr
&+2\sum_{k=1}^\infty\sum_{j=1}^{\vartheta(k)}\alpha_{j,k}t_{j,k}(m)
t_{j,k}(n)(\varphi_+)^+((\txt{1\over2}-k)i;u,v,w,z;g),
}
$$
where
$$
(\varphi_+)^+(r;u,v,w,z;g)
:={{\pi{i}}\over{2\sinh(\pi{r})}}\int_0^\infty
(J_{2ir}(x)-J_{-2ir}(x))\varphi_+(x;u,v,w,z;g){{\d
x}\over{x}},\leqno(3.19)
$$
and $J_\nu(x)$ is the Bessel function of the first kind in standard
notation
(see [L]). Also,
$$
\leqalignno{
K_{-}&(m,n;u,v,w,z;g):=\sum_{j=1}^{\infty}\varepsilon_{j}
\alpha_{j}t_{j}(m)t_{j}(n)
(\varphi_{-})^-(\kappa_{j};u,v,w,z;g)&(3.20)\cr
&+{1\over\pi}\int_{-\infty}^{\infty}{{\sigma_{2ir}(m)\sigma_{2ir}
(n)}\over{(mn)^{ir}|\zeta(1+2ir)|^{2}}}(\varphi_{-})^-
(r;u,v,w,z;g)\d r,
}
$$
where
$$
(\varphi_{-})^-(r;u,v,w,z;g):=2\cosh(\pi r)
\int_{0}^{\infty}\varphi_{-}(r;u,v,w,z;g)K_{2ir}(x){{\d
x}\over{x}},\leqno(3.21)
$$
and $K_\nu(x)$ is the Bessel function of imaginary argument
(or Macdonald's function).

\medskip
Now, in order to facilitate later discussion, we introduce three
functions ${\Phi}_\pm$ and $\Xi$ of five complex variables:
$$
\leqalignno{
{\Phi}_+(\xi;&\,u,v,w,z;g):=-i(2\pi)^{w-z-2}\cos(\hf\pi(u-v))&(3.22)\cr
&\times\int_{-i\infty}^{i\infty}\sin(\hf\pi(u+v+w+z-2s))\cr
&\times\Gamma(\hf(u+v+w+z-1)+\xi-s)
\Gamma(\hf(u+v+w+z-1)-\xi-s)\cr
&\times\Gamma(s+1-u-w)\Gamma(s+1-v-w)\tilde{g}(s,w)\d s;
}
$$
$$
\leqalignno{
\Phi_-(\xi;&\,u,v,w,z;g)=i(2\pi)^{w-z-2}\cos(\pi\xi)
\int_{-i\infty}^{i\infty}\cos(\pi(w+\hf(u+v)-s))&(3.23)\cr
&\times\Gamma(\hf(u+v+w+z-1)+\xi-s)
\Gamma(\hf(u+v+w+z-1)-\xi-s)\cr
&\times\Gamma(s+1-u-w)\Gamma(s+1-v-w)\tilde{g}(s,w)\d s;
}
$$
$$
\leqalignno{
\Xi(\xi;u,v,w,z;g)=&
{1\over{2\pi{i}}}\int_{-\infty i}^{\infty i}
{{\Gamma\left(\xi+{1\over2}(u+v+w+z-1)-s\right)}\over{\Gamma\left(\xi+
{1\over 2}(3-u-v-w-z)+s\right)}}&(3.24)\cr
&\times\Gamma(s+1-u-w)
\Gamma(s+1-v-w)\tilde{g}(s,w)\d s.
}
$$
Note that the path in $(3.22)$ is such that
the poles of the first two gamma-factors and those of the
other three factors in the integrand are separated to the right and the
left, respectively, by the path, and $\xi,u,v,w,z$ are assumed to be
such
that the path can be drawn. The path in $(3.23)$ is chosen in just
the same way.  On the other hand the path in $(3.24)$ separates the
poles
of $\Gamma\left(\xi+{1\over 2}(u+v+w+z-1)-s\right)$ and
those of $\Gamma(s+1-u-w)\Gamma(s+1-v-w)\tilde{g}(s,w)$ to
the left and the right of the path, respectively.

\medskip

{\bf Lemma 4.}
{\it We have
$$
\leqalignno{
\Phi_{+}(\xi;u,v,&w,z;g)=-{{(2\pi)^{w-z}\cos\big({1\over
2}\pi(u-v))}\over{4\sin(\pi\xi)}}&(3.25)\cr
&\times\left\{\Xi(\xi;u,v,w,z;g)-\Xi(-\xi;u,v,w,z;g)\right\};
}
$$
$$
\leqalignno{
\Phi_{-}(\xi;u,v,w,z;g)=&{(2\pi)^{w-z}\over{4\sin(\pi\xi)}}
\Bigg\{\sin(\pi(\hf(z-w)+\xi))\Xi(\xi;u,v,w,z;g)&(3.26)\cr
&-\sin(\pi(\hf(z-w)-\xi))\Xi(-\xi;u,v,w,z;g)\Bigg\},
}
$$
provided the left sides are well-defined. Also, for real $r$
and} $(u,v,w,z)\in{\cal E}$ ({\it see} (3.12)),
$$
(\varphi_+)^+(r;u,v,w,z;g)= \hf(2\pi)^{1-w+z}
\Phi_{+}(ir;u,v,w,z;g),
\leqno(3.27)
$$
$$
(\varphi_{-})^-(r;u,v,w,z;g)= \hf(2\pi)^{1-w+z}
\Phi_{-}(ir;u,v,w,z;g);
\leqno(3.28)
$$
{\it and, for integral $k\ge1$ and} $(u,v,w,z)\in {\cal E}$,
$$
\leqalignno{
(\varphi_{+})^+&(i(\hf-k);u,v,w,z;g)&(3.29)\cr
&= \hf(-1)^{k}\pi\cos(\hf\pi(u-v))\Xi(k-\hf;u,v,w,z;g).
}
$$
The last three formulas are consequences of Mellin transforms of $J$-
and $K$-
Bessel functions.
\medskip
Next, we insert the spectral expansions $(3.18)$ and
$(3.20)$ into $(3.15)$ and exchange the order of sums and integrals.
The absolute convergence that we have to check is obvious as far as the
double
summation over the variables $m,n$ is concerned, since we have
$(u,v,w,z)\in\cal E$ and
$$
t_j(n)\ll n^{{1\over4}+\delta},\quad t_{j,k}(n)\ll
n^{{1\over4}+\delta},\leqno(3.30)
$$
where the implicit constant depends only on $\delta$, an arbitrary
fixed
positive constant. The bounds in (3.30) are not the best ones known,
but they are
sufficient for our purpose. Thus the issue is reduced to bounding
$(\varphi_\pm)^\pm$; and Lemma 4 renders it in terms of the function
$\Xi$. We
then have, uniformly for any fixed compact subset of $\cal{E}$,
$$
\Xi(ir;u,v,w,z;g)\ll |r|^{-{1\over4}A},\quad
\Xi(k-\txt{1\over2};u,v,w,z;g)\ll k^{-{1\over4}A},
\leqno(3.31)
$$
as real $r$ and positive integral $k$ tend to infinity.  Hence, on
noting (see [Mo6]) that
$$
\sum_{K\le\kappa_j<2K}\alpha_j\ll K^2,\quad
\sum_{j=1}^{\vartheta(k)}\alpha_{j,k}
\ll k,\leqno(3.32)
$$
we are now able to  exchange freely the order of sums and integrals in
question,  as long as we work inside $\cal{E}$.
\par
Before stating our new expressions for ${\e{J}}_3^\pm$ we put
$$
\leqalignno{
S(\xi;u&,v,w,z):=\zeta(\hf(u+v+w+z-1)+\xi)\zeta(\hf(u+v+w+z-1)-
\xi)&(3.33)\cr
\times&\zeta(\hf(u+z-v-w+1)+\xi)
\zeta(\hf(u+z-v-w+1)-\xi)\cr\times&\zeta(\hf
(v+z-u-w+1)+\xi)\zeta(\hf(v+z-u-w+1)-\xi).\cr
}
$$
Then we have
\medskip
{\bf Lemma 5.}
{\it In the domain $\cal E$ we have
$$
{\e J}_3^+(u,v,w,z;g)={\e J}_{3,c}^+(u,v,w,z;g)+
{\e J}_{3,d}^+(u,v,w,z;g)+
{\e J}_{3,h}^+(u,v,w,z;g),
\leqno(3.34)
$$
where
$$
{\e J}_{3,c}^+(u,v,w,z;g):={1\over{i\pi}}\int_{(0)}
{{S(\xi;u,v,w,z)}\over{\zeta(1+2\xi)\zeta(1-2\xi)}}
\Phi_+(\xi;u,v,w,z;g)\d\xi,
\leqno(3.35)
$$
$$
\leqalignno{\quad
{\e J}_{3,d}^+(u,v,&w,z;g):=\sum_{j=1}^\infty
\alpha_jH_j(\hf(u+v+w+z-1))H_j(\hf(u+z-v-w+1))&(3.36)\cr
&\times  H_j(\hf(v+z-u-w+1))\Phi_+(i\kappa_j; u,v,w,z;g),
\cr
}
$$
$$
\leqalignno{\quad
&{\e
J}_{3,h}^+(u,v,w,z;g):=(2\pi)^{w-z}\cos(\hf(u-v))&(3.37)\cr
&\times\sum_{k=6}^\infty\sum_{j=1}^{\vartheta(k)}(-1)^k
\alpha_{j,k}H_{j,k}(\hf(u+v+w+z-1))H_{j,k}(\hf(u+z-v-w+1))\cr
&\times H_{j,k}(\hf(v+z-u-w+1))\Xi(k-\hf; u,v,w,z;g).\cr
}
$$
Also
$$
{\e J}_3^-(u,v,w,z;g)=
{\e J}_{3,c}^-(u,v,w,z;g)+{\e J}_{3,d}^-(u,v,w,z;g),
\leqno(3.38)
$$
where
$$
{\e J}_{3,c}^-(u,v,w,z;g):={1\over{i\pi}}\int_{(0)}
{{S(\xi;u,v,w,z)}\over{\zeta(1+2\xi)\zeta(1-2\xi)}}
\Phi_-(\xi;u,v,w,z;g)\d\xi,
\leqno(3.39)
$$
$$
\leqalignno{\qquad
{\e J}_{3,d}^-(u,v,w,z;g)=&\sum_{j=1}^\infty
\varepsilon_j\alpha_jH_j(\hf(u+v+w+z-1))H_j(\hf(u+z-v-w+1))&(3.40)\cr
&\times  H_j(\hf(v+z-u-w+1))\Phi_-(i\kappa_j; u,v,w,z;g).
\cr
}
$$
}
\medskip
Our next task is to show that the above spectral
expansions of ${\e J}_3^\pm$ can be continued to the domain $\cal B$,
whereby
we shall finish our meromorphic continuation of ${\e J}_1$. The domain
$\cal
B$ is obviously symmetric and wide enough to have a joint domain with
${\cal D}_+$, where the decomposition $(3.9)$ was introduced. Hence
$(3.9)$ should hold throughout $\cal B$, and we shall obtain a spectral
decomposition of ${\e L}_4(g;\tau,\tau)$, as asserted.
\medskip
By virtue of Lemma 4, our problem is equivalent to
studying the analytic properties of the function $\Xi\,$. In fact,
it is meromorphic in a fairly wide domain in ${\Bbb C}^5$, as given by
\medskip
{\bf Lemma 6.}
{\it The function $\Xi(\xi;u,v,w,z;g)$ is
meromorphic in the domain
$$
\widetilde{\cal B}=\left\{\xi\,: \Re\xi>-{\txt{1\over8}}A\right\}
\times{\cal{B}}
\leqno(3.41)
$$
and regular in $\widetilde{\cal B}\setminus{\cal{N}}$, where
$\cal N$ is the set of points $(\xi, u,v,w,z)$ such that at least one
of
$$
\xi+\hf(u+v+w+z-1),\quad
\xi+\hf(u+z-v-w+1),\quad
\xi+\hf(v+z-u-w+1)\leqno(3.42)
$$
is equal to a non-positive integer.
Moreover, if $|\xi|$ tends to infinity in any fixed vertical or
horizontal strips while satisfying $\Re\xi>-{1\over8}A$, then
uniformly in ${\cal B}$ we have}
$$
\Xi(\xi;u,v,w,z;g)\ll |\xi|^{-{1\over4}A}.
\leqno(3.43)
$$

In passing, we record that we have also (this is [Mo6, Lemma 4.8]
with $\gamma = 2\xi+1$)
\medskip
{\bf Lemma 7.}
{\it If $(\xi,u,v,w,z)$ is such that the path in $(3.24)$ can be drawn
in a vertical strip contained in the half plane $\Re{s}>0$, then we
have
$$
\leqalignno{
&\Xi(\xi;u,v,w,z;g)
&(3.44)\cr
&={{\Gamma(\alpha)
\Gamma(\beta)}
\over{\Gamma(2\xi+1)}}
\int_0^\infty{{{y}^{\xi+{1\over2}(u+v+w+z
-3)}}\over{(1+y)^{w}}}g^*(\log(1+y))
F(\alpha,\beta;2\xi+1;-y)\d y,
}
$$
where $F$ }({\it see} (2.9)) {\it is the hypergeometric function, and
$$
\alpha = \xi+\hf(u+z-v-w+1),\quad\beta=\xi+\hf(v+z-u-w+1).
\leqno(3.45)
$$
}
\medskip

An immediate consequence of Lemma 6 is that ${\e{J}}_{3,d}^\pm$ and
$\e{J}_{3,h}^+$ are meromorphic inside $\cal B$. Thus, we shall
consider ${\e
J}_{3,c}^\pm$.  To this end we assume first that $(u,v,w,z)$ is in
${\cal E}$; and put
$$
{\e J}_{3,c}(u,v,w,z;g)={\e J}_{3,c}^+(u,v,w,z;g)
+{\e J}_{3,c}^-(u,v,w,z;g).
\leqno(3.46)
$$
We have, by $(3.25)$--$(3.26)$,
$$
\leqalignno{
&{\e J}_{3,c}(u,v,w,z;g)=i(2\pi)^{w-z-1}\int_{(0)}
{{S(\xi;u,v,w,z)}\over{\sin(\pi\xi)\zeta(1+2\xi)\zeta(1
-2\xi)}}&(3.47)\cr
\times&\{\cos(\hf\pi(u-v))
-\sin(\pi(\hf(z-w)+\xi))\}\Xi(\xi;u,v,w,z;g)\d\xi.\cr
}
$$
Applying the functional equation for $\z(s)$ to $\zeta(1-2\xi)$, we
obtain from (3.47)
$$
\leqalignno{\qquad
&{\e J}_{3,c}(u,v,w,z;g)&(3.48)\cr
&=2i(2\pi)^{w-z-2}\int_{(0)}(2\pi)^{2\xi}
\{\cos(\hf\pi(u-v))
-\sin(\pi(\hf(z-w)+\xi))\}\cr&\times
S(\xi;u,v,w,z)\Gamma(1-2\xi)\{\zeta(2\xi)
\zeta(1+2\xi)\}^{-1}\Xi(\xi;u,v,w,z;g)\d\xi.\cr
}
$$
We then choose  $Q$ which is to satisfy the condition
$$
3B<Q\le{\txt{1\over4}A};\quad\zeta(s)
\ne0\quad\hbox{\rm for $\Im{s}=\pm{Q}$}.\leqno(3.49)
$$
We divide the range of integration in $(3.48)$ into two parts according
as $|\xi|>Q$ and $|\xi|\le{Q}$, and denote the corresponding parts of
${\e J}_{3,c}$ by ${\e J}_{3,c}^{(1)}$ and ${\e J}_{3,c}^{(2)}$,
respectively. We observe that if $\Re\xi=0$, $|\Im\xi|\ge{Q}$, then
$S(\xi;u,v,w,z)$ is regular and $O(|\xi|^{cB})$ uniformly in
${\cal B}$ with an absolute constant $c$. Then,
Lemma 7 implies that the integrand in the part ${\e J}_{3,c}^{(1)}$ is
regular
and of fast decay with respect to $\xi$ uniformly in ${\cal B}$.  Hence
${\e J}_{3,c}^{(1)}$ is regular in
${\cal B}$.  As to ${\e{J}}_{3,c}^{(2)}$, we move the path to $L_Q$
which
is the result of connecting the points $-Qi$, $[Q]+{1\over4}-Qi$,
$[Q]+{1\over4}+Qi$, $Qi$ with straight lines. The singularities of the
integrand which we encounter in this procedure are all poles, and
located at
$$
\hf(u+v+w+z-3),\quad \hf(u+z-v-w-1),\quad\hf(v+z-u-w-1);
\leqno(3.50)
$$
$$
\hf\rho\quad(\,|\Im\rho|<Q\,); \quad \hf n\quad
(\,2\le{n}\le[Q]\,);
\leqno(3.51)
$$
where $\rho$ is a complex zero of $\zeta(s)$; note our choice of $Q$.
The first three come from
$S(\xi;u,v,w,z)$, and the others from $\Gamma(1-2\xi)\zeta(2\xi)^{-1}$,
since we have here
$(u,v,w,z)\in {\cal E}$ and so the $\Xi$-factor is regular for
$\Re\xi\ge0$.
We may suppose, for an obvious reason, that the poles given in
$(3.50)$ are all simple, and do not coincide with any of those given in
$(3.51)$. Then we have
$$
{\e J}_{3,c}(u,v,w,z;g)=F_-(u,v,w,z;g)+U(u,v,w,z;g)+
{\e J}_{3,c}^{(Q)}(u,v,w,z;g).
\leqno(3.52)
$$
Here $F_-$ and $U$ are the contributions of residues at the poles given
in $(3.50)$ and $(3.51)$, respectively; and ${\e J}_{3,c}^{(Q)}$ is the
same as $(3.48)$ but with the path $L_Q^*$ which is the sum of the
path $L_Q$ and the half lines $(-i\infty,-Qi]$, $[Qi,i\infty)$. By
virtue of
Lemma 7, the terms $F_-$ and $U$ are meromorphic over $\cal B$, and
${\e
{J}}_{3,c}^{(Q)}$ is regular there.
\medskip
Summing up, we have
\medskip
{\bf Lemma 8.}
{\it The function} ${\e J}_1(u,v,w,z;g) $ (see (3.11))
{\it continues meromorphically to the domain $\cal{B}$. Thus the
decomposition
$(3.9)$ holds throughout $\cal B$.
}
\bigskip
It remains for us
only to specialize $(3.9)$ by setting $(u,v,w,z)=\r{p}_\tau$. This
amounts to studying the local behaviour, near $\r{p}_\tau$, of the
various components of ${\e J}_1$ which have been introduced in the
above
discussion. As a consequence we shall obtain the explicit
formula for ${\e{L}}_4(g;\tau,\tau)$ furnished by Theorem 1.
Namely we have the decomposition, over $\cal B$,
$$
{\e J}_1(u,v,w,z;g)=\{{\e J}_2+{\e J}_{3,c}+{\e J}_{3,d}^+
+{\e J}_{3,d}^-+{\e J}_{3,h}^+\}(u,v,w,z;g),
\leqno(3.53)
$$
where  ${\e{J}}_{3,d}^\pm$
is  regular at $\r{p}_\tau$. To see this we observe that when
$(u,v,w,z)$
is near $\r{p}_\tau$ the point $(ir,u,v,w,z)$, with an arbitrary real
$r$,
is  not in the set ${\cal{N}}$ defined at $(3.42)$; thus by
$(3.25)$--$(3.26)$ the
functions $\Phi_\pm(ir;u,v,w,z)$ are also regular at $\r{p}_\tau$ for
any real
$r$. Hence ${\e{J}}_{3,d}^\pm$ are regular at $\r{p}_\tau$. Similarly
one can
see that ${\e{J}}_{3,h}^+$ is regular at $\r{p}_\tau$. That is, we may
set
$(u,v,w,z)=\r{p}_\tau$ in the series expansions (3.36), (3.37) and
(3.40) without any modification, and find that
$$
\leqalignno{
\{{\e J}_{3,d}^+&+{\e J}_{3,d}^-+{\e
J}_{3,h}^+\}(\r{p}_\tau;g)&(3.54)\cr
&=\sum_{j=1}^\infty\alpha_jH_j^2(\hf)H_j(2\tau-\hf)
\{\Phi_++\Phi_-\}(i\kappa_j;\r{p}_\tau;g)\cr
&+\sum_{k=1}^\infty\sum_{j=1}^{\vartheta(2k)}
\alpha_{j,2k}H_{j,2k}^2(\hf)H_{j,2k}(2\tau-\hf)
\Xi(2k-\hf;\r{p}_\tau;g),\cr
}
$$
where we have used the fact that
$H_j({1\over2})=0$ if $\varepsilon_j=-1$, and
$H_{j,k}({1\over2})=0$ if $k$ is odd.
Note also that by $(3.25)$--$(3.26)$ we have, for real $r$,
$$
\leqalignno{
\{\Phi_++\Phi_-\}(ir;\r{p}_\tau;g)=&
{1\over4}\left(1+{i\over{\sinh(\pi{r})}}\right)\Xi(ir;\r{p}_\tau;
g)&(3.55)\cr
+&{1\over4}\left(1-{i\over{\sinh(\pi{r})}}\right)\Xi(-ir;\r{p}_\tau;g).
}
$$
On the other hand, Lemma 7 gives
$$
\leqalignno{
\Xi(ir;\r{p}_\tau;g)=&{{\Gamma({1\over2}+ir)^2}\over{\Gamma(1+2ir)}}
\int_0^\infty
y^{2\tau-{3\over2}+ir}(1+y)^{-\tau}{g}^*(\log(1+y))&(3.56)\cr
&\times{F}(\hf+ir,\hf+ir;1+2ir;-y)\d y.\cr
}
$$
Hence
$$
\leqalignno{
&\{\Phi_++\Phi_-\}(ir;\r{p}_\tau;g)={1\over2}\int_0^\infty
y^{2\tau-{3\over2}}(1+y)^{-\tau}g^*(\log(1+y))&(3.57)\cr
&\times{\Re}\left\{y^{ir}\Big(1+{i\over{\sinh(\pi{r})}}\Big)
{{\Gamma({1\over2}+ir)^2}\over{
\Gamma(1+2ir)}}F(\hf+ir,\hf+ir;1+2ir;-y)
\right\}\d y.\cr
}
$$
Further, we observe that $(3.56)$ holds with $ir$ replaced by
$k-{1\over2}$; and thus (3.57) gives, for any integer $k\ge1$,
$$
\{\Phi_++\Phi_-\}(2k- \hf;\r{p}_\tau;g)
=\Xi(2k-\hf;\r{p}_\tau;g).
\leqno(3.58)
$$
\medskip
Now, we consider ${\e J}_{3,c}$ in an immediate neighbourhood of
$\r{p}_\tau$. Let us assume first that
$$
{1\over2}<\tau<{3\over4}.\leqno(3.59)
$$
\smallskip
\noindent
This is much similar to the case $\tau={1\over2}$, which is treated in
[Mo6]. We
return to $(3.52)$, and move the contour in
${\e {J}}_{3,c}^{(Q)}$ back to the imaginary axis, while keeping
$(u,v,w,z)$
close to $\r{p}_\tau$. The poles which we encounter in this process are
those given in $(3.51)$ and ${1\over2}(3-u-v-w-z)$, which is in fact to
the right
of the imaginary axis. Other poles of $S(\xi;u,v,w,z)$ are either
on the left of the imaginary axis or cancelled by the zeros of the
factor
$\cos({1\over2}\pi(u-v))-\sin(\pi({1\over2}(z-w)+\xi))$, and moreover
Lemma 7
implies that
$\Xi(\xi;u,v,w,z;g)$ is regular for ${\Re}(\xi)\ge-{1\over4}$. We
denote by $F_+(u,v,w,z;g)$ the contribution of the pole
${1\over2}(3-u-v-w-z)$. Then we have
$$
{\e{J}}_{3,c}^{(Q)}(u,v,w,z;g)=F_+(u,v,w,z;g)
-U(u,v,w,z;g)+{\e J}_{3,c}^*(u,v,w,z;g),\leqno(3.60)
$$
where ${\e J}_{3,c}^*$ has the same expression as the right side of
$(3.48)$ but with different $(u,v,w,z)$. Hence, by $(3.52)$,
$$
{\e J}_{3,c}(u,v,w,z;g)=\{F_++F_-\}(u,v,w,z;g)+{\e
J}_{3,c}^*(u,v,w,z;g),
\leqno(3.61)
$$
when $(u,v,w,z)$ is close to $\r{p}_\tau$. Here we should note that
${\e J}_{3,c}^*$ is regular at $\r{p}_\tau$, and
$$
{\e J}_{3,c}^*(\r{p}_\tau;g)={1\over\pi}\int_{-\infty}^\infty
{{|\zeta({1\over2}+it)|^4|\zeta(2\tau-
{1\over2}+it)|^2}\over{|\zeta(1+2it)|^2}}
\{\Phi_++\Phi_-\}(it;\r{p}_\tau;g)\d t.
\leqno(3.62)
$$
This ends the local study of the decomposition $(3.53)$ in the vicinity
of
$\r{p}_\tau$, provided that $(3.59)$ holds.
\medskip
Now, if  $(u,v,w,z)$ is close to $\r{p}_\tau$, then we have
$$
\leqalignno{
{\e J}&(u,v,w,z;g)=M(u,v,w,z;g)+{\e J}_{3,c}^*(u,v,w,z;g)+
\overline{{\e{J}}_{3,c}^*(\overline{w},\overline{z},
\overline{u},\overline{v};g)}&(3.63)\cr
&+\{{\e{J}}_{3,d}^-+{\e J}_{3,d}^++{\e
J}_{3,h}^+\}(u,v,w,z;g)+\overline{
\{{\e{J}}_{3,d}^-+{\e J}_{3,d}^++{\e
J}_{3,h}^+\}(\overline{w},\overline{z},\overline{u},
\overline{v};g)}\,,
}
$$
where
$$
\leqalignno{
M&(u,v,w,z;g)={\e{J}}_0(u,v,w,z;g)+{\e{J}}_2(u,v,w,z;g)+
\overline{{\e{J}}_2(\overline{w},\overline{z},\overline{u},
\overline{v};g)}&(3.64)\cr
&+\{F_++F_-\}(u,v,w,z;g)+\overline{
\{F_++F_-\}(\overline{w},\overline{z},\overline{u},
\overline{v};g)}.
}
$$
It should be stressed that all terms in $(3.63)$ are regular at
$\r{p}_\tau$. That the function $M$ is regular at $\r{p}_\tau$ is due
to the fact that all terms in $(3.63)$ except for $M$ have already
been proved to be regular at $\r{p}_\tau$.
\par
It remains for us to express
$M(\r{p}_\tau;g)$ in terms of $g$. We have
$$
M(u,v,w,z;g)\;=\;\sum_{j=0}^{12}M_j(u,v,w,z;g)\leqno(3.65)
$$
with
$$
M_{6+j}(u,v,w,z;g)=\overline{M_j
(\overline{w},\overline{z},\overline{u},
\overline{v};g)}\qquad (\,1\le{j}\le6\,)\,.\leqno(3.66)
$$
Here (recall that $\tilde g$ is given by (3.1))
$$
M_0(u,v,w,z;g)=g^*(0)\zeta(u+w)\zeta(u+z)\zeta(v+w)\zeta(v+z)
\{\zeta(u+v+w+z)\}^{-1},\leqno(3.67)
$$
$$
\leqalignno{
M_1(u,v,w,z;g)&=\tilde{g}(v+w-1,w)\zeta(u+z)\zeta(v+w-1)&(3.68)\cr
&\times\zeta(z-w+1)\zeta(u-v+1)\{\zeta(u+z-v-w+2)\}^{-1},
\cr}
$$
$$\leqalignno{
M_2(u,v,w,z;g)&=\tilde{g}(u+w-1,w)\zeta(v+z)\zeta(u+w-1)&(3.69)\cr
&\times\zeta(z-w+1)\zeta(v-u+1)\{\zeta(v+z-u-w+2)\}^{-1},
\cr
}
$$
$$
\leqalignno{
M_3(u,v,w,z;g)&=(2\pi)^{w-z}\{\cos(\hf\pi(u-v))
+\cos(\pi (z-w+ \hf
(u-v)))\}&(3.70)\cr
&\times\zeta(u+z-1)\zeta(v+w)
\zeta(z-w)\zeta(v-u+1)\cr
&\times\{\cos(\hf
\pi(u+z-v-w))\zeta(2-u-z+v+w)\}^{-1}\cr
&\times\Xi(\hf(u+z-v-w-1);u,v,w,z;g),\cr
}
$$
$$
\leqalignno{
M_4(u,v,w,z;g)&=(2\pi)^{w-z}\{\cos(\hf\pi(u-v))
+\cos(\pi(z-w+\hf(v-u)))\}&(3.71)\cr
&\times\zeta(v+z-1)\zeta(u+w)
\zeta(z-w)\zeta(u-v+1)\cr
&\times\{\cos(\hf\pi(v+z-u-w))\zeta(2-v-z+u+w)\}^{-1}\cr
&\times\Xi(\hf(v+z-u-w-1);u,v,w,z;g),\cr
}
$$
$$
\leqalignno{
M_5(u,v,w,z;g)&=-(2\pi)^{w-z}\{\cos(\hf\pi(u-v))
-\cos(\pi(z+\hf(u+v)))\}&(3.72)\cr
&\times\zeta(u+z-1)\zeta(2-v-w)
\zeta(v+z-1)\zeta(2-u-w)\cr
&\times\{\cos(\hf\pi(u+v+w+z))\zeta(4-u-v-w-z)\}^{-1}\cr
&\times\Xi(\hf(u+v+w+z-3);u,v,w,z;g),\cr
}
$$
$$
\leqalignno{
M_6(u,v,w,z;g)&=(2\pi)^{w-z}\{\cos(\hf\pi(u-v))
+\cos(\pi(w+\hf(u+v)))\}&(3.73)\cr
&\times\zeta(u+z-1)\zeta(2-v-w)
\zeta(v+z-1)\zeta(2-u-w)\cr
&\times\{\cos(\hf
\pi(u+v+w+z))\zeta(4-u-v-w-z)\}^{-1}\cr
&\times\Xi(-\hf(u+v+w+z-3);u,v,w,z;g)\,.\cr
}
$$
Among these, $M_0$ is equal to ${\e J}_0$; $M_1$ and
$M_2$ come from ${\e {J}}_2$; and $M_j$ $(\,3\le{j}\le6\,)$
are the contributions of residues of the integral in $(3.48)$ at the
poles
$\xi={1\over2}(u+z-v-w-1)$, ${1\over2}(v+z-u-w-1)$,
${1\over2}(u+v+w+z-3)$,
${1\over2}(3-u-v-w-z)$, respectively. They can be singular at
$\r{p}_\tau$
individually, but the singular parts should cancel each other out if
they are
brought into $(3.64)$, for $M$ is regular at $\r{p}_\tau$. More
precisely,
put $(u,v,w,z)=\r{p}_\tau+(a_1,a_2,a_3,a_4)\delta$ with a small
complex $\delta$, and expand each term into a Laurent series in
$\delta$;
then the sum of the constant terms is equal to
$M(\r{p}_\tau)$, regardless of the choice of the vector
$(a_1,a_2,a_3,a_4)$. We choose it in such a way that it is real and no
singularities of any of the $M_j$ $(0\le j\le12)$ are encountered
when $|\delta|$ tends to $0$. This is possible, for the exceptional
$a_1,a_2,a_3,a_4$ satisfy a finite number of linear relations. Thus we
shall
assume hereafter that $\delta\ne0$ is small and the vector
$(a_1,a_2,a_3,a_4)$ is chosen accordingly; and we denote
$(a_1,a_2,a_3,a_4)\delta$ either by $(\delta)$ or by
$(\delta_1,\delta_2,\delta_3,\delta_4)$. Also we denote the constant
term
of $M_j$ by $M_j^*$.
\medskip
First, we have trivially
$$
M_0^*=M_0(\r{p}_\tau;g) \;=\;
{\zeta^4(2\tau)\over\zeta(4\tau)}g^*(0).\leqno(3.74)
$$
Invoking $(3.1)$, we have
$$
\leqalignno{
M_1(\r{p}_\tau+(\delta);g)&=\Gamma(2\tau-1+\delta_2+\delta_3)
\zeta(2\tau+\delta_1+\delta_4)\zeta(2\tau-1+\delta_2+\delta_3)
&(3.75)\cr
&\times \zeta(\delta_4-\delta_3+1)\zeta(\delta_1-\delta_2+1)
\{\zeta(2+\delta_1-\delta_2-\delta_3+\delta_4)\}^{-1}\cr
&\times\int_{-\infty}^\infty
{{\Gamma(1-\tau-\delta_2-it)}\over\Gamma(\tau+\delta_3-it)}
g(t)\d t.\cr
}
$$
This implies that the singularity of $M_1$ at $\r{p}_\tau$ is
of order two. Hence the constant term of
$M_1(\r{p}_\tau+(\delta);g)$ is a linear combination of the first three
coefficients of the power series in $\delta$ for the last integral.
Thus
$$
\leqalignno{
&\qquad
M_1^*=\int_{-\infty}^\infty{\Gamma(1-\tau-it)\over\Gamma(\tau-it)}
&(3.76)\cr
&\times \left(d_0+d_1{\Gamma'
\over\Gamma}(1-\tau-it){\Gamma'\over\Gamma}
(\tau-it)+d_2{\Gamma''
\over\Gamma}(1-\tau-it)+d_3{\Gamma''
\over\Gamma}(\tau-it)\right)g(t)\d t
}
$$
where the constants $d_j$ depend on $\tau$ and the vector
$(a_1,a_2,a_3,a_4)$. Clearly $M_2$ can be treated
in just the same way, and $M_2^*$ has the same form as $(3.76)$.
\par
The terms $M_j$ $(\,3\le{j}\le6\,)$ are not so simple; and our
computation of them depends on a classical formula of Barnes (see e.g.,
[WW]).
By the definition $(3.24)$ we have, for the $\Xi$-factor in $M_3$,
$$
\leqalignno{
&\Xi(\hf(u+z-v-w-1);u,v,w,z;g)&(3.77)\cr
&={1\over{2\pi{i}}}
\int_{-i\infty}^{i\infty}\Gamma(u+z-1-s)\Gamma(s+1-u-w)\tilde{g}(s,w)\d
s,
\cr}
$$
where the path separates the poles of $\Gamma(u+z-1-s)$ and those of
the other
two factors to the right and the left, respectively; that we can draw
such
a path is assured by our choice of $(a_1,a_2,a_3,a_4)$. Inserting
$(3.1)$ in this we get an absolutely convergent double integral,
hence it follows that
$$
\leqalignno{
\Xi&(\hf(u+z-v-w-1);u,v,w,z;g)=
\int_{-\infty}^\infty{{g(t)}\over{\Gamma(w-it)}}&(3.78)\cr
&\times{1\over{2\pi{i}}}\int_{-i\infty}^{i\infty}
\Gamma(s)\Gamma(s+1-u-w)\Gamma(u+z-1-s)
\Gamma(w-it-s)\d s\d t.\cr
}
$$
The path of the inner integral is the same as in $(3.1)$; and obviously
we
may suppose that it separates the poles of the first two
$\Gamma$-factors
from those of the other two. Hence we have, again by the  Barnes
formula,
$$
\leqalignno{
\Xi&(\hf(u+z-v-w-1);u,v,w,z;g)&(3.79)\cr
&=\Gamma(u+z-1)\Gamma(z-w)
\int_{-\infty}^\infty{{\Gamma(1-u-it)}\over{\Gamma(z-it)}}g(t)\d t.
}
$$
This implies that $M_3$ has a singularity of order two at $\r{p}_\tau$;
thus
$M_3^*$ admits an expression of the same form as $(3.76)$.
Obviously the same argument applies to $M_4$.
\par
The
$\Xi$-factor of $M_5$ can be computed in much the same way, and we have
$$
\leqalignno{
\Xi&(\hf(u+v+w+z-3);u,v,w,z;g)&(3.80)\cr
&=\Gamma(u+z-1)\Gamma(v+z-1)
\int_{-\infty}^\infty{{\Gamma(1-u-it)\Gamma(1-v-it)}\over{\Gamma(w-it)
\Gamma(z-it)}}g(t)\d t.\cr
}
$$
This implies that $M_5$ is regular at $\r{p}_\tau$, and
$$
\leqalignno{
M^*_5=M_5(\r{p}_\tau;g)&=(1-\sec(2\tau\pi)){(\Gamma(2\tau
-1)\zeta(2\tau-1)
\zeta(2-2\tau))^2\over\zeta(4-4\tau)}&(3.81)\cr
&\times\int_{-\infty}^\infty\left({\Gamma(1-\tau-it)\over\Gamma(\tau-
it)}
\right)^2g(t)\d t
}
$$
A rearrangement gives
$$
\leqalignno{
&\int_{-\infty}^\infty\left({\Gamma(1-\tau-it)\over\Gamma(\tau-it)}
\right)^2g(t)\d t&(3.82)\cr
&={1\over2\pi^2}\int_{-\infty}^\infty
|\Gamma(1-\tau+it)|^4(1-\cos(2\tau\pi)\cosh(2\pi t))g(t)\d t.
}
$$
\par
As to $M_6$, this also is regular at $\r{p}_\tau$, since
the $\Xi$-factor is regular there because of Lemma 7. We have
$$
M_6^*=M_6(\r{p}_\tau;g)=(1+\sec(2\tau\pi)){(\zeta(2\tau-1)
\zeta(2-2\tau))^2\over\zeta(4-4\tau)}\Xi({\txt{3\over2}}-2\tau;
\r{p}_\tau;g)\leqno(3.83)
$$
with
$$
\Xi({\txt{3\over2}}-2\tau;\r{p}_\tau;g)=
{1\over2i}\int_{-\infty}^\infty{g(t)\over\Gamma(\tau-it)}
\left(\int_{(\sigma_0)}{\Gamma(s+1-2\tau)^2\Gamma(\tau-it-
s)\over\sin(\pi s)
\Gamma(3-4\tau+s)}\d s\right)\d t,\leqno(3.84)
$$
where $2\tau-1<\sigma_0<\tau$. We have
$$
\leqalignno{
\int_{(\sigma_0)}\cdots=-{2\pi i\over\sin(\pi(\tau-it))}&\sum_{j=1}^N
{\Gamma(j+1-2\tau)^2\over\Gamma(j+3-4\tau)
\Gamma(j+1-\tau+it)}&(3.85)\cr
+{2\pi i\over\sin(\pi(\tau-it))}&\sum_{j=0}^{N-1}{\Gamma(j+1-\tau-it)^2
\over\Gamma(j+3(1-\tau)-it)\Gamma(j+1)}+\int_{(\sigma_N)}\cdots.
}
$$
where $2\tau-1+N<\sigma_N<\tau+N$.
This ends the discussion under the assumption
that $(3.59)$ holds.

\medskip
Next, let
$$
\tau={3\over4}.\leqno(3.85)
$$
There is an essential difference between this case and $(3.59)$,
which we just discussed.
This is due to the fact that  the singularities
$\pm{1\over2}(u+v+w+z-3)$ of
the integrand in $(3.47)$ approach  the origin as $(u,v,w,z)$ tends
to
$\r{p}_{3\over4}$. That is, they cannot be treated as well-separated.
\par
While keeping $(u,v,w,z)$ close to $\r{p}_{3\over4}$, we move the
contour in
${\e {J}}_{3,c}^{(Q)}$ to the imaginary axis but with a small outward
indent
around the origin. The poles which we encounter in this
process are those given in $(3.51)$ and $\pm{1\over2}(u+v+w+z-3)$.
Other poles
of $S(\xi;u,v,w,z)$ are either close to
$-{1\over2}$ or cancelled by the zeros of the factor
$\cos({1\over2}\pi(u-v))-
\sin(\pi({1\over2}(z-w)+\xi))$, and moreover Lemma 7 implies that
$\Xi(\xi;u,v,w,z;g)$ is regular for ${\Re}(\xi)\ge-{1\over4}$. We have
$$
{\e{J}}_{3,c}^{(Q)}(u,v,w,z;g)=(F_+-F_-)(u,v,w,z;g)
-U(u,v,w,z;g)+{\e J}_{3,c}^*(u,v,w,z;g),\leqno(3.86)
$$
where $F_\pm$ are as before, and ${\e J}_{3,c}^*$ has the same
expression as the
right side of $(3.47)$ but with the indented contour and a different
$(u,v,w,z)$.
By $(3.52)$,
$$
{\e J}_{3,c}(u,v,w,z;g)=F_+(u,v,w,z;g)+{\e J}_{3,c}^*(u,v,w,z;g),
$$
when $(u,v,w,z)$ is close to $\r{p}_{3\over4}$. Here
we should note that
${\e J}_{3,c}^*$ is regular at $\r{p}_{3\over4}$, and
$$
{\e J}_{3,c}^*(\r{p}_{3\over4};g)={1\over\pi}\int_{-\infty}^\infty
{{|\zeta({1\over2}+it)|^4|\zeta(1+it)|^2}\over{|\zeta(1+2it)|^2}}
\{\Phi_++\Phi_-\}(it;\r{p}_{3\over4};g)\d t,
\leqno(3.87)
$$
because this integrand is continuous.
This ends the local study of the decomposition $(3.53)$ in a small
neighbourhood
of $\r{p}_{3\over4}$.
\medskip
Now, if  $(u,v,w,z)$ is close to $\r{p}_{3\over4}$, then the
counterpart of
$(3.63)$--$(3.64)$ holds, and it remains for us to express
$M(\r{p}_{3\over4};g)$ in terms of $g$, but with the new $M$. We have
$$
M(u,v,w,z;g)=\sum_{\scr{j=0}\atop\scr{j\ne5,11}}^{12}
M_j(u,v,w,z;g)\leqno(3.88)
$$
where $M_j$ are the same as in $(3.65)$--$(3.73)$. The terms $M_5$ and
$M_{11}$
are missing, because the shift of the contour cancels the contribution
of the pole at ${1\over2}(u+v+w+z-3)$ out, as we have seen
above.
\par
The computation of $M_j^*$ is the same as before. It should perhaps be
remarked that
$$
M_6^*=M_6(\r{p}_{3\over4};g)=-{2\over\pi}\zeta^4(\hf)
\Xi\left(0;\r{p}_{3\over4};g\right).\leqno(3.89)
$$
We have
$$
\leqalignno{
\qquad\Xi\left(0;\r{p}_{3\over4};g\right)&=
{1\over2\pi i}\int_{-\infty}^\infty{g(t)\over\Gamma({3\over4}-it)}
\left(\int_{({2\over3})}\Gamma^2(s-\hf)\Gamma({\txt{3\over4}}-it-s)
\Gamma(1-s)\d s\right)\d t&(3.90)\cr
&=\pi\int_{-\infty}^\infty\left({\Gamma({1\over4}-it)
\over\Gamma({3\over4}-it)}\right)^2
g(t)\d t.
}
$$
\medskip
Finally, let
$$
{3\over4}<\tau<1.\leqno(3.91)
$$
Then the pole ${1\over2}(3-u-v-w-z)$ is on the left of the imaginary
axis; and
the contribution of the pole ${1\over2}(u+v+w+z-3)$ is cancelled out by
moving the contour to the imaginary axis. That is, we have
$$
M(u,v,w,z;g)=\sum_{\scr{j=0}\atop\scr{j\ne5,6,11,12}}^{12}
M_j(u,v,w,z;g)\leqno(3.92)
$$
This ends our discussion and completes the proof of Theorem 1.

\bigskip
\centerline{\bb 4. Sums of spectral values}
\bigskip
Note that (2.6) of Theorem 1 contains the quantities
$\a_j H_j^2(\hf)H_j(\s)\;(\s = 2\tau - \hf)$
with a given  $\hf < \s < 1$, while in Part II of this work we shall
encounter sums containing $\a_j H_j(\hf)H_j^2(\s)$. For the
omega-results relating to moments we shall need the non-vanishing
of
$$
L_\s(\k) \;:=\;\sum_{\k_j=\k}\a_j H^2_j(\hf)H_j(\s),\quad
N_\s(\k) \;:=\;\sum_{\k_j=\k}\a_j H_j(\hf)H_j^2(\s)\leqno(4.1)
$$
for infinitely many $\kappa$ and a given  $\hf < \s < 1$. The
non-vanishing of $L_{1\over2}(\k)$ was used (see [I2], [I5], [I7],
[IM1], [Mo4],
[Mo6]) for omega
results on the fourth moment of $|\zt|$. The non-vanishing
of $L_\s(\k)$ and $N_\s(\k)$  that we
need is a corollary of the following

\medskip
THEOREM 2. {\it For fixed $\tau$ such that $\hf < \tau < 1$
and $K\to\infty$, we have
$$
\sum_{\kappa_j\le K}\alpha_jH_j(\hf)H_j^2(\tau)=
(1+o(1))\pi^{-2}\zeta^2\!(\tau+\hf)
\zeta(2\tau)K^2\leqno(4.2)
$$
and}
$$
\sum_{\kappa_j\le K}\alpha_jH_j^2(\hf)H_j(\tau)=
(1+o(1))2\pi^{-2}\zeta^2\!(\tau+\hf)
K^2\log K.\leqno(4.3)
$$

\medskip
For $L_{1\over2}(\k)$ not only that non-vanishing is known,
but a sharper asymptotic formula for the sum in
question, namely
$$
\sum_{\k_j\le K} \a_jH_j^3(\hf) =
K^2P_3(\log K) + O(K^{5/4}\log^{37/4}K),\leqno(4.4)
$$
proved by the first author [Iv9], where $P_3(x)$ is a suitable
cubic polynomial.
One could also employ similar methods to obtain a sharpening of (4.2)
and (4.3),
but this will not be done here, since it is not needed in
the sequel. It is known (see Katok--Sarnak [KS]) that $H_j(\hf)
\ge 0$; it follows trivially from the functional equation for
$H_j(s)$ that $H_j(\hf)=0$ if $\eps_j = -1$.
Our formula (4.3) supports the conjecture that
$H_j(\s)\ge0$ for $0\le\s\le1$, but
this remains an open problem.

\medskip
{\bf Proof of Theorem 2}.
Because $H_j\left({1\over2}\right)=0$ when $\eps_j=-1$,
we may start by treating the sum
$$
\sum_{\kappa_j<K}\left(\eps_j\alpha_j
H_j(\hf)H_j(\lambda)\right)\cdot H_j(\tau),
\leqno(4.5)
$$
with the aim of taking later $\lambda={1\over2}$ or $\lambda = \tau$
in (4.5).
\medskip
Let $h(r)$ be an even, entire function such that $h(\pm\hf i) = 0$
and $h(r) \ll\exp(-c|r|^2)\;(c>0)$ in any fixed horizontal strip,
and
$$
H(u,v;f;h):=\sum_{j=1}^\infty\eps_j\alpha_jH_j(u)H_j(v)t_j(f)h(\kappa_j)
\qquad(f\ge1).\leqno(4.6) $$
Transformation formulas for the sums appearing in (4.6) were
established by the
second author [Mo1] and then in [Mo6] (see eq.\ (3.3.6) there, also
$(3.3.8)$ and $(3.3.9)$ are important). The formulas in question
transform the quantity (4.6) from spectral theory into a sum of
various quantities from classical analytic number theory. We set
$$ \e{H}(\lambda;f;h)=H(\hf,\lambda;f;h).\leqno(4.7) $$
According to the formulas displayed on p.\ 117 of [Mo6], we have
$$
\leqalignno{ &H(u,v;f;h)=
2(\pi{i})^{-1}(2\pi\sqrt{f})^{2(u-1)}{{\hat{h}(1-u)}\over{\cos(\pi{u})}}
\sigma_{1-u-v}(f)\zeta(1-u+v)&(4.8)\cr
&+2(\pi{i})^{-1}(2\pi\sqrt{f})^{2(v-1)}{{\hat{h}(1-
v)}\over{\cos(\pi{v})}}
\sigma_{1-u-v}(f)\zeta(1-v+u)\cr
&+8(2\pi)^{u+v-4}\sum_{m=1}^\infty{m}^{u-1}\sigma_{v-u}(m)\sigma_{1-u-
v}(m+f)
\Psi_+(u,v;m/f;h)\cr
&+8(2\pi)^{u+v-4}\sum_{\scr{m=1}\atop\scr{m\ne{f}}}^\infty
{m}^{u-1}\sigma_{v-u}(m)\sigma_{1-u-v}(m-f) \Psi_-(u,v;m/f;h)\cr
&+8(2\pi)^{u+v-4}{f}^{u-1}\sigma_{v-u}(f)\zeta(u+v-1)
\Psi_-(u,v;1;h)\cr
&-4\sigma_{2(u-1)}(f)f^{1-u}\zeta(u+v-1)\zeta(v-u+1)h(i(u-1))/\zeta(3
-2u)\cr
&-4\sigma_{2(v-1)}(f)f^{1-v}\zeta(u+v-1)\zeta(u-v+1)h(i(v-1))/\zeta(3
-2v)\cr
&-\pi^{-1}\int_{-\infty}^\infty{{\sigma_{2ir}(f)
\zeta(u+ir)\zeta(u-ir)\zeta(v+ir)
\zeta(v-ir)}\over{f^{ir}|\zeta(1+2ir)|^{2}}}h(r)\d r.\cr } $$ Here
$$ \leqalignno {
\Psi_+(u,v;x;h)=-&\int_{(\beta)}\Gamma(1-u-s)\Gamma(1-v-s)&(4.9)\cr
&\times
\cos(\pi(s+\txt{1\over2}(u+v))){{\hat{h}(s)}\over{\cos(\pi{s})}}x^s\d
s } $$ and $$ \leqalignno{ \Psi_-&(u,v;x;h)
=\cos(\txt{1\over2}\pi(u-v))&(4.10)\cr &\times\int_{(\beta)}
\Gamma(1-u-s)\Gamma(1-v-s) {{\hat{h}(s)}\over{\cos(\pi{s})}}x^s\d
s, } $$ where $$ \hat{h}(s)=\int_{\Im
r=-C}rh(r){\Gamma(s+ir)\over\Gamma(1-s+ir)}\d r\qquad (\Re
s>-C),\leqno(4.11)
$$
with any large $C>0$. This is proved if
$u\ne{v}$, which can be dropped by an obvious convention; and
$1+\beta<\Re(u),\Re(v)<1$ with $-{3\over2}<\beta<0$.
\par
If
$$
{1\over2}<\lambda<1,\leqno(4.12)
$$
then we have
$$
\leqalignno{
\lim_{(u,v)\to\left({1\over2},\lambda\right)}&\Big\{(2\pi\sqrt{f})^{2(u-
1)}
{{\hat{h}(1-u)}\over{\cos(\pi{u})}}
\sigma_{1-u-v}(f)\zeta(1-u+v)&(4.13)\cr
&+(2\pi\sqrt{f})^{2(v-1)}{{\hat{h}(1-v)}\over{\cos(\pi{v})}}
\sigma_{1-u-v}(f)\zeta(1-v+u)\Big\}\cr
=&{1\over2\pi^2}\zeta(\lambda+\hf)(\hat{h})'(\hf)
\sigma_{{1\over2}-\lambda}(f)f^{-{1\over2}}\cr
+&(2\pi)^{2(\lambda-1)}\sec(\pi\lambda)\zeta({\txt{3\over2}}-\lambda)
\hat{h}(1-\lambda)\sigma_{{1\over2}-\lambda}(f)f^{\lambda-1},
}
$$
where the fact $\hat{h}\left({1\over2}\right)=0$ has been used
(see $(3.3.15)$ of  [Mo6]), and that
$$
\leqalignno{
\lim_{(u,v)\to\left({1\over2},\lambda\right)}&\big\{\sigma_{2(u
-1)}(f)f^{1-u}
\zeta(u+v-1)\zeta(v-u+1)h(i(u-1))/\zeta(3-2u)&(4.14)\cr
&+\sigma_{2(v-1)}(f)f^{1-v}\zeta(u+v-1)
\zeta(u-v+1)h(i(v-1))/\zeta(3-2v)\big\}\cr
=&{6\over\pi^2}\zeta(\lambda-\hf)\zeta(\lambda+\hf)h-\hf i)
\sigma_{-1}(f)f^{1\over2}\cr
&+{\zeta\left(\lambda-{1\over2}\right)
\zeta\left({3\over2}-\lambda\right)\over\zeta(3-2\lambda)}h(i(\lambda
-1))
\sigma_{2(\lambda-1)}(f)f^{1-\lambda}.
}
$$
Note that the right sides of both $(4.13)$ and $(4.14)$ have removable
singularities
at $\lambda={1\over2}$. It is found that if $(4.12)$ holds, then
$$
\e{H}(\lambda;f;h)=\sum_{\nu=1}^7\e{H}_\nu(\lambda;f;h).\leqno(4.15)
$$
Here
$$
\e{H}_1(\lambda;f;h)={2\over\pi i}\times\hbox{the right side of
$(4.13)$},
\leqno(4.16)
$$
$$
{\e{H}}_2(\lambda;f;h)=8(2\pi)^{\lambda-{7\over2}}
\sum_{m=1}^\infty{m}^{-{1\over2}}
\sigma_{\lambda-{1\over2}}(m)\sigma_{{1\over2}-\lambda}(m+f)
\Psi^+(\lambda;m/f;h),\leqno(4.17)
$$
$$
{\e{H}}_3(\lambda;f;h)=8(2\pi)^{\lambda-{7\over2}}
\sum_{m=1}^\infty(m+f)^{-{1\over2}}
\sigma_{\lambda-{1\over2}}(m+f)
\sigma_{{1\over2}-\lambda}(m)\Psi^-(\lambda;1+m/f;h),\leqno(4.18)
$$
$$
{\e{H}}_4(\lambda;f;h)=8(2\pi)^{\lambda-{7\over2}}\sum_{m=1}^{f-1}m^{-
{1\over2}}
\sigma_{\lambda-{1\over2}}(m)\sigma_{{1\over2}-\lambda}(f-m)
\Psi^-(\lambda;m/f;h),\leqno(4.19)
$$
$$
{\e{H}}_5(\lambda;f;h)=8(2\pi)^{\lambda-{7\over2}}
f^{-{1\over2}}\sigma_{{1\over2}-\lambda}(f)
\zeta(\lambda-\hf)\Psi^-(\lambda;1;h),
\leqno(4.20)
$$
$$
{\e{H}}_6(\lambda;f;h)=-4\times\hbox{the right side of
$(4.14)$},\leqno(4.21)
$$
$$
{\e{H}}_7(\lambda;f;h)=-\pi^{-1}\int_{-
\infty}^\infty{{|\zeta({1\over2}+ir)|^2
|\zeta(\lambda+ir)|^2}\over
{|\zeta(1+2ir)|^2}}\sigma_{2ir}(f)f^{-ir}h(r)\d r,\leqno(4.22)
$$
where
$$
\Psi^+(\lambda;x;h)=-\int_{(\beta)}
\Gamma(\hf-s)\Gamma(1-\lambda-s)
{\cos(\pi(s+{1\over2}(\lambda+{1\over2}))\over
\cos\pi s}\hat{h}(s)x^s\d s,
\leqno(4.23)
$$
$$
\Psi^-(x;h)=\cos(\hf\pi(\hf-\lambda))
\int_{(\beta)}\Gamma(\hf-s)\Gamma(1-\lambda-
s){{\hat{h}(s)}\over{\cos\pi{s}}}x^s\d s
\leqno(4.24)
$$
with $-{3\over2}<\beta<1-\lambda$. If we consider the limit as
$\lambda$ tends to
${1\over2}$, then from (4.15)--(4.24) we obtain the
assertion of Lemma 3.8 of [Mo6].
\medskip
We have
$$
(\hat{h})'(\hf)=2\int_{-\infty}^\infty
rh(r){\Gamma'\over\Gamma}(\hf+ir)\d r,\leqno(4.25)
$$
and
$$
\hat{h}(1-\lambda)=\int_{-\infty}^\infty
rh(r){\Gamma(1-\lambda+ir)\over\Gamma(\lambda+ir)}\d r.\leqno(4.26)
$$
Also,
$$
\leqalignno{
&\qquad\Psi^+(\lambda;x;h)=-{1\over\pi i}\int_{-\infty}^\infty
rh(r)\sinh\pi
r&(4.27)\cr &\times\int_{(\beta)}\Gamma(\hf-s)\Gamma(1-\lambda-s)
\Gamma(s+ir)\Gamma(s-ir)
\cos\{\pi(s+\hf(\lambda+\hf))\}x^s\d s \d r,
}
$$
with $0<\beta<1-\lambda$. Evaluating the inner integral, we have, for
any $x>0$,
$$
\leqalignno{
&\Psi^+(\lambda;x;h)=-2\pi i\int_{-\infty}^\infty {rh(r)\over \cosh\pi
r}
\cos\pi(ir-\hf(\lambda+\hf))&(4.28)\cr
&\times{\Gamma({1\over2}+ir)\Gamma(1-\lambda+ir)\over\Gamma(1+2ir)}
F(\hf+ir,1-\lambda+ir;1+2ir;-1/x)
x^{-ir} \d r.
}
$$
Then, by Gauss' integral representation for $F$ (see e.g., [L]
or [WW]),
$$
\leqalignno{
&\Psi^+(\lambda;x;h)=-2\pi i\int_0^1\left\{y(1-y)\right\}^{-{1\over2}}
\left(1+{y\over x}\right)^{\lambda-1}&(4.29)\cr
&\times\int_{-\infty}^\infty rh(r)
{\cos\pi\left(ir-
{1\over2}\left(\lambda+{1\over2}\right)\right)\over\cosh\pi r}
{\Gamma(1-\lambda+ir)\over\Gamma({1\over2}+ir)}\left({y(1-y)\over
x+y}\right)^{ir}\d r \d y,
}
$$
which corresponds to $(3.3.41)$ of [Mo6].
\par
In what concerns $\Psi^-$, for $x>1$ we obtain in a similar
fashion
$$
\leqalignno{
&\Psi^-(\lambda;x;h)=2\pi
i\cos(\hf\pi(\hf-\lambda))\int_{-\infty}^\infty{rh(r)\over
\cosh(\pi r)}&(4.30)\cr
&\times{\Gamma({1\over2}+ir)\Gamma(1-\lambda+ir)\over\Gamma(1+2ir)}
F(\hf+ir,1-\lambda+ir;1+2ir;1/x)x^{-ir} \d r\cr
&=2i\cos(\hf\pi(\hf-\lambda))\int_0^1\left\{y(1-y)\right\}^{-{1\over2}}
\left(1-{y\over x}\right)^{\lambda-1}\cr
&\times\int_{-\infty}^\infty rh(r)\Gamma(\hf-ir)
\Gamma(1-\lambda+ir)\left({y(1-y)\over
x-y}\right)^{ir}\d r \d y,
}
$$
which corresponds to $(3.3.43)$ of [Mo6]. Also,
$$
\Psi^-(\lambda;1;h)=
{2\pi\cos(\hf\pi(\hf-\lambda))
\over\Gamma({3\over2}-\lambda)}\int_{-\infty}^\infty
rh(r)\tanh(\pi r)\Gamma(1-\lambda+ir)\Gamma(1-\lambda-ir)\d
r.\leqno(4.31)
$$
When $0<x<1$, we argue as on p.\ 121 of [Mo6], to deduce that
$$
\leqalignno{
\Psi^-(x;h)&=\cos(\hf\pi(\hf-\lambda))\int_0^\infty
\left\{\int_{-\infty}^\infty{r}h(r)\left({y\over{1+y}}\right)^{ir}\d
r\right\}
&(4.32)\cr&\times \left\{\int_{(\beta)}x^s(y(y+1))^{s-1}
{{\Gamma({1\over2}-s)\Gamma(1-\lambda-s)}
\over{\Gamma(1-2s)\cos(\pi{s})}}\d s\right\}\d y,
}
$$
with $-{3\over2}<\beta<1-\lambda$, $\beta\ne-{1\over2}$.
\bigskip
We shall now derive  an approximate functional equation for
$H_j(\tau)$. The expression
$(4.5)$ implies in particular that we may restrict ourselves to the
case
$\eps_j=+1$. Let us assume that
$$
|\kappa_j-K|\le{G}\log{K},\quad(\log{K})^2<G<K^{1-\delta}\qquad
(\delta>0).\leqno(4.33)
$$
Take a large $C>0$ and consider the integral
$$
{\cal{R}}:={1\over{2\pi{i}\mu}}\int_{(3)}H_j(w+\tau)K^w
\Gamma(w/\mu)\d w\qquad (\mu=C\log K).
\leqno(4.34)
$$
We have
$$
{\cal{R}}=\sum_{f\le3K}t_j(f)f^{-\tau}\exp(-(f/K)^\mu)+O({\rm
e}^{-K}).\leqno(4.35)
$$
Shifting the path of integration in (4.34) to
$\Re w = -\hf\mu$ and
recalling the functional equation for $H_j(s)$, we obtain
$$
{\cal{R}}=H_j(\tau)+\sum_{f=1}^\infty{t_j}(f)f^{\tau-1}
{\cal{R}}_j(fK),\leqno(4.36)
$$
where
$$
\leqalignno{
{\cal{R}}_j(x):={1\over{(2\pi)^{2(1-\tau)}}\pi i\mu}
\int\limits_{\left(-{1\over2}\mu\right)}&
(4\pi^2x)^w\Gamma(1-\tau-w+i\kappa_j)
\Gamma(1-\tau-w-i\kappa_j)&(4.37)\cr
&\times\big\{\cosh(\pi\kappa_j)
-\cos(\pi(w+\tau))\big\}\Gamma(w/\mu)\d w.
}
$$
By Stirling's formula for the gamma-function the above integrand
is
$$
\ll(4\pi^2x)^{-{1\over2}\mu}
\big(|w+i\kappa_j||w-i\kappa_j|\big)^{{1\over2}
(\mu+1)-\tau}\exp(-\pi|w|/(2\mu))\leqno(4.38)
$$
and thus
$$
{\cal{R}}_j(x)=O\left(K^{1-2\tau}(4\pi^2xK^{-2})^{-{1\over2}\mu}
\right),\leqno(4.39)
$$
where the implied constant is absolute. This allows us to truncate the
last sum
over $f$ at $f=[3K]$ with an error which is $\ll K^{-C}$ for any
fixed $C>0$. Hence we have proved

\medskip
{\bf Lemma 10}. {\it For $\hf < \tau < 1$ fixed and
  uniformly for all $\kappa_j$ satisfying $(4.33)$ and $\eps_j = 1$,
$$
\leqalignno{
H_j(\tau)=&\sum_{f\le{3K}}t_j(f)f^{-\tau}\exp(-(f/K)^\mu)&(4.40)\cr
&-\sum_{f\le{3K}}t_j(f)f^{\tau-1}
{\cal{R}}_j^{(1)}(fK)+O(K^{-{1\over2}C}),\cr
}
$$
where $C>0$ is any given constant and
$$
\leqalignno{
\quad{\cal{R}}_j^{(1)}(x):={1\over{(2\pi)^{2(1-\tau)}}\pi
i\mu}&\int\limits_{-\mu^{-1}-i\mu^2} ^{-\mu^{-1}+i\mu^2}
(4\pi^2x)^w\Gamma(1-\tau-w+i\kappa_j)\Gamma(1-\tau-w-i\kappa_j)
&(4.41)\cr
&\times\big\{\cosh(\pi\kappa_j)
-\cos(\pi(w+\tau))\big\}\Gamma(w/\lambda)\d w.\cr
}
$$}

\medskip
Stirling's formula gives, for any $N\ge1$
and for the values of $w$ relevant in $(4.41)$,
$$
\leqalignno{
\log&\,\Gamma(1-\tau-w+i\kappa_j)
=(\hf-\tau-w+i\kappa_j)\log(1-\tau-w+i\kappa_j)&(4.42)\cr
&+\tau-1+w-i\kappa_j
+\txt{1\over2}\log(2\pi)
+\sum_{\nu=1}^{2N}b_\nu(1-\tau-w+i\kappa_j)^{-\nu}
+O(K^{-2N-{1\over2}}),\cr
}
$$
where $b_\nu$'s are absolute constants, and the implied constant
depends only
on $N$. Therefore
$$
\leqalignno{
\log&\,\Gamma(1-\tau-w+i\kappa_j)=
(\hf-\tau-w+i\kappa_j)\big\{\log(\kappa_j)
+\txt{1\over2}\pi{i}\big\}&(4.43)\cr
&-i\kappa_j+\txt{1\over2}\log(2\pi)
+\sum_{\nu=1}^{2N}{p_\nu}(w)\kappa_j^{-\nu}+O(K^{
-2N}(\log{K})^{12N+6})\cr
}
$$
with certain polynomials $p_\nu$ of degree $\le\nu+1$ with constant
coefficients. Adding to this the corresponding formula for
$\Gamma(1-\tau-w-i\kappa_j)$, we have
$$
\leqalignno{
\log\big\{\Gamma(1-\tau-w+i\kappa_j)&
\Gamma(1-\tau-w-i\kappa_j)\big\}=(1-2\tau-2w)\log(\kappa_j)
-\pi\kappa_j&(4.44)\cr&+\log(2\pi)
+\sum_{\nu=1}^N{p_{2\nu}}(w)\kappa_j^{-2\nu}+O(K^{
-2N}(\log{K})^{12N+6}).\cr
}
$$
This implies readily that the integrand of $(4.41)$ can be replaced by
$$
\pi\kappa_j^{1-2\tau}
(4\pi^{2}\kappa_j^{-2}x)^w\Big\{1+\sum_{\nu=1}^N{q_\nu}(w)\kappa_j^{
-2\nu}
+O(K^{-N})\Big\}\Gamma(w/\lambda),\leqno(4.45)
$$
where $q_\nu(w)$ are polynomials of degree $\le3\nu$ with constant
coefficients,
and the $O$-constant depends only on $N$. Then we expand each
$\kappa_j^{1-2\tau-2w-2\nu}$ into a power series in
$(1-(\kappa_j/K)^2)=O(K^{-\delta}\log{K})$ and truncate it at the power
$N_1=[2N/\delta]$. Rearranging the result of truncation we
see that the integrand of $(4.41)$ can be written as
$$
\pi K^{1-2\tau}(4\pi^2K^{-2}x)^w\big\{Q(w,1-(\kappa_j/K)^2)+O(K^{-N})
\big\}\Gamma(w/\lambda),
\leqno(4.46)
$$
where
$$
Q(w,y)=\sum_{\nu=0}^{N_1}u_\nu(w)y^\nu,\;
u_0(w)=1+\sum_{\nu=1}^N{q}_\nu(w){K}^{-2\nu}.\leqno(4.47)
$$
Inserting
$(4.47)$ into $(4.41)$ and restoring the range of integration to the
whole line
$\Re{w}=-\mu^{-1}$, we get, uniformly for $f\le{3K}$,
$$
{\cal{R}}_j^{(1)}(fK)=\left({K\over\pi}\right)^{1-2\tau}
\sum_{\nu=0}^{N_1}U_\nu(fK)
(1-(\kappa_j/K)^2)^\nu+O(K^{-N}).\leqno(4.48)
$$
Here $N_1=[2N/\delta]$ and
$$
U_\nu(x)={1\over{2\pi{i}\mu}}\int_{(-\mu^{-1})}(4\pi^2K^{-2}x)^wu_\nu(w)
\Gamma(w/\mu)\d w,
\leqno(4.49)
$$
where $u_p(w)$ is a polynomial of degree $\le2N_1$, whose coefficients
are
independent of $\kappa_j$ and bounded by a constant depending only
on $\delta$, $\tau$, and $N$. Hence, if $\eps_j=1$, we have, for any
$N\ge1$ and $\mu=C\log{K}$ with a sufficiently large $C>0$,
$$
\leqalignno{
&H_j(\tau)=\sum_{f\le3K}t_j(f)f^{-\tau}\exp(-(f/K)^\mu)&(4.50)\cr
&-\left({K\over\pi}\right)^{1-2\tau}
\sum_{\nu=0}^{N_1}\sum_{f\le{3K}}t_j(f)f^{\tau-1}
U_\nu(fK)(1-(\kappa_j/K)^2)^\nu+O(K^{-{1\over5}N}+K^{-{1\over2}C})
}
$$
with the implied constant depending only on $\delta$, $\tau$, $C$, and
$N$.
\bigskip
We shall evaluate asymptotically, as $K\to\infty$,
$$
\e{C}(\lambda,\tau;K,G) \;:=\; \sum_{j=1}^\infty\alpha_jH_j(\hf)
H_j(\lambda)H_j(\tau)h_0(\kappa_j),\leqno(4.51)
$$
where ${1\over2}<\lambda,\,\tau<1$ initially, and the weight
function will be
$$
h_0(r):=\left(r^2+{1\over4}\right)\left\{
\exp\left(-\left({r-K\over G}\right)^2\right)+
\exp\left(-\left({r+K\over G}\right)^2\right)\right\},\leqno(4.52)
$$
provided that $(4.33)$ holds. From (4.50) we obtain
$$
\leqalignno{
&{\e{C}}(\lambda,\tau;K,G)=
\sum_{f\le3K}f^{-\tau}\exp(-(f/K)^\mu){\e{H}}(\lambda;f;h_0)&(4.53)\cr
&-\left({K\over\pi}\right)^{1-2\tau}
\sum_{\nu=0}^{N_1}\sum_{f\le{3K}}f^{\tau-1}
U_\nu(fK){\e{H}}(\lambda;f;h_\nu)+o(1),
}
$$
where $\e{H}$ is defined by $(4.7)$, and
$$
h_{\nu}(r)=h_0(r)(1-(r/K)^2)^\nu.\leqno(4.54)
$$
\par
To evaluate $\e{H}(\lambda;f;h_\nu)$, we use $(4.15)$. The
contributions of $(4.17)$,
$(4.18)$, $(4.20)$ and $(4.21)$ are negligible, which can be confirmed
in much the
same way as on pp.\ 128--129 of [Mo6]. Then, corresponding to
$(3.4.25)$ there, we
have ($d(n) \equiv \s_0(n)$ is the number of divisors of $n$)
$$
\leqalignno{
{\e{H}}&(\lambda;f;h_\nu)=\e{H}_1(\lambda;f;h_\nu)
+O\left(d(f)(G+K^{2\over3})K^2(G/K)^\nu(\log K)^c\right), &(4.55)\cr
+&O\left(f^{3\over2}K^3G^{-3}(G/K)^\nu\log
K\sum_{m<f}m^{-2}\sigma_{\lambda-{1\over2}}(m)
\sigma_{{1\over2}-\lambda}(f-m)\right)}
$$
with some constant $c>0$, provided that
$$
K^{{1\over2}+\delta}<G<K^{1-\delta}.\leqno(4.56)
$$
Note that $(4.55)$ holds with $\lambda={1\over2}$, too, and that the
estimation of the error terms is not the best that our argument can
attain.
Also, we should remark that
$$
\e{H}_1(\lambda;f;h_\nu)\ll K^3G(G/K)^\nu(\log fK)
\sigma_{{1\over2}-\lambda}(f)f^{-{1\over2}},\leqno(4.57)
$$
uniformly for ${1\over2}\le\lambda<1$ and for all $f\ge1$.
In fact, when ${1\over2}+(\log K)^{-1}\le
\lambda<1$, this follows from $(4.16)$, $(4.25)$, $(4.26)$, and
otherwise one may use
the Taylor expansion at $\lambda={1\over2}$.
\par
Inserting $(4.55)$ in $(4.53)$ we obtain
$$
\leqalignno{
\qquad\e{C}(\lambda,\tau;K,G)&=\sum_{f<3K}\e{H}_1(\lambda;f;h_0)
\left(f^{-\tau}\exp(-(f/K)^\mu)-\left({K\over\pi}\right)^{1
-2\tau}f^{\tau-1}
U_0(fK)\right)&(4.58)\cr
&+O\left(K^{3-\tau}G(\log K)^c\right)\cr
&=\sum_{f=1}^\infty\e{H}_1(\lambda;f;h_0)
\left(f^{-\tau}\exp(-(f/K)^\mu)-\left({K\over\pi}\right)^{1
-2\tau}f^{\tau-1}
U_0(fK)\right)\cr
&+O\left(K^{3-\tau}G(\log K)^c\right),
}
$$
provided that ${1\over2}\le\lambda,\,\tau<1$ and
$$
K^{{2\over3}+\delta}<G<K^{1-\delta}.\leqno(4.59)
$$
The extension of the summation to $f\ge 3K$ can be performed in view of
$(4.49)$
with an appropriate shift of the contour to the left.
\par
This means that we have
$$
\leqalignno{
&\sum_{f=1}^\infty\e{H}_1(\lambda;f;h_0)f^{-\tau}\exp(-(f/
K)^\mu)&(4.60)\cr
=&-{1\over2\pi^4\mu}\zeta(\lambda+\hf)(\hat{h}_0)'(\hf)\cr
&\qquad\times\int_{(3)}\zeta(w+\tau+\hf)
\zeta(w+\lambda+\tau)K^w\Gamma\left({w\over\mu}\right)\d w\cr
&-{1\over\pi^2\mu}(2\pi)^{2(\lambda-1)}\sec(\pi\lambda)
\zeta({\txt{3\over2}}-\lambda)\hat{h}_0(1-\lambda)\cr
&\qquad\times\int_{(3)}\zeta(w+\tau+\hf)
\zeta(w+1-\lambda+\tau)K^w\Gamma\left({w\over\mu}\right)\d w.
}
$$
Also,
$$
\leqalignno{
&\sum_{f=1}^\infty\e{H}_1(\lambda;f;h_0)f^{\tau-1}U_0(fK)&(4.61)\cr
=&-{1\over2\pi^4\mu}\zeta(\lambda+\hf)(\hat{h}_0)'(\hf)\cr
&\qquad\times\int_{(-3)}\zeta({\txt{3\over2}}-\tau-w)
\zeta(1+\lambda-\tau-w)(4\pi^2/
K)^wu_0(w)\Gamma\left({w\over\mu}\right)\d w\cr
&-{1\over\pi^2\mu}(2\pi)^{2(\lambda-1)}\sec(\pi\lambda)
\zeta({\txt{3\over2}}-\lambda)
\hat{h}_0(1-\lambda)\cr
&\qquad\times\int_{(-3)}\zeta({\txt{3\over2}}-\tau-w)
\zeta(2-\lambda-\tau-w)(4\pi^2/
K)^wu_0(w)\Gamma\left({w\over\mu}\right)\d w.
}
$$
\bigskip
Before specialising the above formula, note that
$$
(\hat{h}_0)'(\hf)=2i\pi^{3\over2}K^3G+O(KG^3).
\leqno(4.62)
$$
We then put $\lambda=\tau$, ${1\over2}<\tau<1$. Then the right side of
$(4.60)$ is
asymptotically equal to
$$
{2\over\pi^{3\over2}}\zeta^2\!(\tau+\hf)\zeta(2\tau)
K^3G,\leqno(4.63)
$$
and that of $(4.61)$ to
$$
-{2\over\pi^{3\over2}}\zeta^2\!(\tau+\hf)
K^3G(4\pi^2/K)^{{1\over2}-\tau}u_0(\hf-\tau)\Gamma
\Bigl({\hf-\tau\over\mu}\Bigr).\leqno(4.64)
$$
Inserting these expressions into $(4.60)$ we find that
$$
\e{C}(\tau,\tau;K,G)=(1+o(1)){2\over\pi^{3\over2}}
\zeta^2\!(\tau+\hf)\zeta(2\tau)K^3G,\leqno(4.65)
$$
which leads to (4.2). Namely, similarly as in [Iv9, eq.
(7.10)-(7.11)], we note that we have
$$\eqalign{&
\int\limits_{K_0}^{2K_0}\e{C}(\tau,\tau;K,G)\d K
= \sum_{j\ge1}\a_j H_j(\hf)H_j^2(\tau)\int_{K_0}^{2K_0}
(\k_j^2+{\txt{1\over4}})\exp(-(\k_j-K)^2G^{-2})\d K + O(1)
\cr&
= \sqrt{\pi}G\sum_{K_0<\k_j \le 2K_0}\a_j H_j(\hf)H_j^2(\tau)\k_j^2
+ o(K_0^4G).\cr}
$$
On the other hand, from the main term on the right-hand side
of (4.65) we obtain
$$
\eqalign{&
2\pi^{-3/2}\z^2(\tau+\hf)\z(2\tau)\int_{K_0}^{2K_0}
K^3G\d K\cr&
= \hf G\pi^{-3/2}\z^2(\tau+\hf)\z(2\tau)((2K_0)^4-K_0^4).\cr}
$$
Here we take $G = K_0^{1-\eps}$ say, then we replace $K_0$ by
$K_02^{-\ell}$ and sum over $\ell\ge1$, and finally replace
$K_0$ by $K$ to obtain
$$
\sum_{\k_j\le K}\a_j H_j(\hf)H_j^2(\tau)\k_j^2 =
\left(\hf\pi^{-2}\z^2(\tau+\hf)\z(2\tau)+o(1)\right)K^4\leqno(4.66)
$$
as $K\to\infty$. The desired formula (4.2) follows then by partial
summation from (4.66).
\medskip
To prove (4.3), set $\lambda={1\over2}$, ${1\over2}<\tau<1$. This
case is treated in Section 3.3 of [Mo6], and we could appeal to
Lemma 3.8 therein. But it is the same as to use $(4.60)$ and
$(4.61)$ with this specialisation. Thus, the right side of
$(4.60)$ with $\lambda={1\over2}$ is equal to
$$ \leqalignno{
{1\over2\pi^4\mu}&\int_{(3)}\left\{2\left(\hat{h}_0\right)'\left({1\over
2}\right) \left(\log(2\pi)-c_E-
{\zeta'\over\zeta}\left(w+\tau+{1\over2}\right)\right)
-{1\over2}\left(\hat{h}_0\right)''\left({1\over2}\right)\right\}&(4.67)
\cr &\times\zeta^2\!\left(w+\tau+{1\over2}\right)
K^w\Gamma\left({w\over\mu}\right)\d w, } $$ and that of $(4.61)$
to
$$ \leqalignno{ {1\over2\pi^4\mu}&\int\limits_{(
-3)}\left\{2\left(\hat{h}_0\right)'\left({1\over2}\right)
\left(\log(2\pi)-c_E-{\zeta'\over\zeta}\left({3\over2}-\tau-
w\right)\right)
-{1\over2}\left(\hat{h}_0\right)''\left({1\over2}\right)\right\}&(4.68)
\cr &\times\zeta^2\!\left({3\over2}-\tau-w\right)
(4\pi^2/K)^wu_0(w)\Gamma\left({w\over\mu}\right)\d w, } $$
 where
$c_E$ is the Euler constant, the $u_0$ is specialized accordingly,
and $$ (\hat{h}_0)''(\hf)=8i\pi^{3\over2}K^3G\log K + O(KG^3\log
K).\leqno(4.69) $$ Hence we have $$
\e{C}(\hf,\tau;K,G)=(1+o(1)){4\over\pi^{3\over2}}
\zeta^2\!(\tau+\hf)K^3G\log K, $$ which implies (4.3) by the
procedure used in the previous case. This completes the proof of
Theorem 2.

\bigskip
\centerline{\bb 5. The asymptotics of the $\Lambda$-function}
\bigskip
We shall apply now Theorem 1 with a specific (Gaussian) exponential
weight
function, namely
$$
g(t) = {1\over2\sqrt{\pi}G}\left\{\exp\left(-\left({T-t\over
G}\right)^2
\right) + \exp\left(-\left({T+t\over G}\right)^2
\right)\right\},\leqno(5.1)
$$
which is a standard one, either in this or in a slightly changed form
(without the factor $1/(2\sqrt{\pi}G)$). Obviously this choice
of $g$ satisfies the {\it basic assumption} in Section 1.

\medskip
The crucial thing needed in the estimation of $E_2(T,\s)$ and related
quantities is the function $\Lambda(r;\tau,g)$, defined by (2.8),
and we proceed in this section to give its asymptotic evaluation.
The main formula is (5.14), but we have found it more
expedient to leave it in this form than to formulate a concrete theorem
or
lemma which would provide the needed asymptotic evaluation.
The form that will be given in the sequel is sharper and more complete
than the one that can be found in [Mo6, Chapter 5]. We suppose  that
the parameters $r,G$ satisfy
$$ 1 \ll r \le
TG^{-1}\log^5T,\quad T^\eps \le G \le T^{1-\eps}, \leqno(5.2)
$$
which are the relevant ranges for our investigations. The case $r<0$
is completely analogous, and the range for $r$ not covered by
(5.2) is treated in [Mo6], where it is shown that the
contribution is negligible. In the  case of the weight function
(5.1) (without the factor $1/(2\sqrt{\pi}G)$) we shall have
$$ g_c(x) =
2\sqrt{\pi}G{\rm e}^{-{1\over4}G^2x^2}\cos(xT).\leqno(5.3)
$$
However, to keep in tune with the notation of [Mo6], we
omit $2\sqrt{\pi}G$ in subsequent calculations. Moreover, the
exponential factor in (5.3) shows that the contribution of $y >
G^{-1}\log T$ in (3.57) is negligible, so that by changing $y$
to $1/y$ it is sufficient to start with the evaluation of the integral
$$
\leqalignno{ I :&= \int_0^{G^{-1}\log
T}y^{2\tau-3/2}(1+y)^{-\tau}\cos(T\log(1+y))\exp(-{\txt{1\over4}}G^2
\log^2(1+y))&(5.4)\cr& \times\Re
\left\{y^{ir}{\G^2(\hf+ir)\over\G(1+2ir)}F(\hf+ir,\hf+ir;
1+2ir;-y)\right\}\d y,\cr}
$$
  where $\tau \ge \hf\;(\tau\ne1)$ is a given constant,
  and of course $I$ depends on $T,r,G$ and $\tau$. There
are several ways to evaluate $I$ asymptotically, but the simplest
procedure seems to use the following quadratic transformation formula
(see
[L, eq. (9.6.12)]), which is valid if
$|\arg(1-z)| < \pi,\, 2\b \not = -1, -3, -5, \ldots\,:$
$$
\leqalignno{ F(\a,\b;2\b;z) &=
{\left({1+\sqrt{1-z}\over2}\right)}^{-2\a}&(5.5)\cr& \times
F\left(\a,\a-\b + \hf; \b+\hf;
{\left({1-\sqrt{1-z}\over1+\sqrt{1-z}}\right)}^2\right).\cr}
$$
Then the relevant part of $I$ becomes
$$
\leqalignno{ &\qquad
\int_0^{G^{-1}\log
T}y^{2\tau-3/2}(1+y)^{-\tau}\cos(T\log(1+y))\exp(-{\txt{1\over4}}G^2
\log^2(1+y))&(5.6)\cr& \times\Re
\left\{y^{ir}{\G^2(\hf+ir)\over\G(1+2ir)}
{\left({1+\sqrt{1+y}\over2}\right)}^{-1-2ir}
F\left(\hf+ir,\hf;
1+ir;{\left({1-\sqrt{1+y}\over1+\sqrt{1+y}}\right)}^2\right)\right\}\d
y.
}
$$
We recall (2.9), and insert it in (5.6) with $\a = \hf + ir,\, \b =
\hf,\, \gamma =
1 + ir$,
$$
z = {\left({1-\sqrt{1+y}\over1+\sqrt{1+y}}\right)}^2
\ll G^{-2}\log^2T = o(1)\quad(T\to\infty),
$$
since  $0 \le y \le G^{-1}\log T$ in (5.6). Note that, for $k \ge 1$,
$$
\left|{(\hf+ir)_k\over(1+ir)_k}\right| \le 1, \quad
{(\hf+ir)_k\over(1+ir)_k} = {\left({ir\over k+ir}\right)}^{1/2}
\left(1 + O\left({1\over r}\right)\right)
$$
uniformly in $k$, with an appropriate choice of branch. Therefore the
main
contribution to $I$ will come from the constant term (i.e., unity) in
the series
expansion (2.9), while the remaining terms will be of a similar nature,
only of a
lower order of magnitude. The series can be truncated in such a way
that the tails
will make a negligible contribution; this procedure will be repeatedly
used
without further explicit mention in subsequent calculations. For
example, we
develop into series the terms $(1+y)^{-\tau}$ and
${\left({1+\sqrt{1+y}\over2}\right)}^{-1}$, noting that the main
contribution will again come from the constant term unity.
Now we use Stirling's formula for the gamma-function
in the form ($t \ge t_0 > 0,\;0\le\s\le1$)
$$
\G(s) = \sqrt{2\pi}\,t^{\s-{1\over2}}\exp\left(-\hf\pi t + it\log t -it
+ \hf{\pi i}(\s - \hf)\right)\cdot\left(1 + O_\s\left(
t^{-1}\right)\right),\leqno(5.7)
$$
with the understanding that the $O$--term in (5.7) admits an
asymptotic expansion in terms of negative powers of $\,t$. Therefore we
have
$$
{\G^2(\hf+ir)\over\G(1+2ir)} =
\sqrt{\pi}r^{-1/2}{\rm e}^{-2ir\log2-{1\over4}\pi i}
\cdot\left(1 + O\left({1\over r}\right)\right)
$$
for the gamma-factors in (5.6),
where  the $O$-term admits an asymptotic expansion. In this
way the problem is reduced to the evaluation of the integral
$$
\leqalignno{ &\sqrt{\pi}r^{-1/2} \int_0^{G^{-1}\log
T}y^{2\tau-3/2}\cos(T\log(1+y))\exp(-{\txt{1\over4}}G^2
\log^2(1+y))&(5.8)\cr& \times\Re\left\{
y^{ir}\exp\left(-2ir\log2-2ir\log
\left({1+\sqrt{1+y}\over2}\right) - {\txt{1\over4}}\pi i\right)\right\}
\d y\cr& =  \sqrt{\pi}r^{-1/2}\int_0^{G^{-1}\log
T}y^{2\tau-3/2}\cos(T\log(1+y))\exp(-{\txt{1\over4}}G^2\log^2(1+y))
\cr& \times
\cos\left(r\log y - r\log4 -
2r\log\left({1+\sqrt{1+y}\over2}\right) - {\txt{1\over4}}\pi \right)
\d y,\cr}
$$
But as
$$
\cos\a\cos\b = \hf[\cos(\a+\b) + \cos(\a-\b)],
$$
we have in fact to consider
$$
\sqrt{\pi}r^{-1/2}{\rm e}^{-ir\log4}
    \int_0^{G^{-1}\log T}y^{2\tau-3/2}
\exp(-{\txt{1\over4}}G^2\log^2(1+y)){\rm
e}^{i{\cal F}_\pm(y,r)-{1\over4}i\pi}\d y,
\leqno(5.9)
$$
with
$$
{\cal F}_\pm(y;r,T) := r\log y - 2r\log\left({1+\sqrt{1+y}\over2}\right)
\pm T\log(1+y),  \leqno(5.10)
$$
so that
$$
{\partial {\cal F}_\pm(y;r,T)\over\partial y} = {r\over y}
  - {r\over1+y+\sqrt{1+y}} \pm {T\over1+y}.
$$
Note that in our range for $y$, which is $0 < y \le G^{-1}\log T$,
the derivative of
${\cal F}_+$ is positive, so there will be no saddle point. Hence we
shall discuss in detail only the more difficult case
of ${\cal F}_-$ (henceforth denoted by ${\cal F}$),
which has a saddle point $y_0$, the root of
$$
{r\over y}  - {r\over1+y+\sqrt{1+y}} =  {T\over1+y}.
$$
This is equivalent to $T^2y^2 - r^2y - r^2 = 0$, giving
$$
y_0 = {r\over T}\left(\sqrt{1 + {r^2\over4T^2}} +
{r\over2T}\right), \leqno(5.11)
$$
so that $y_0 \sim r/T$ as $T\to\infty$. Then
$$
{\cal F}(y_0) = r\log y_0 - 2r\log\left({1+\sqrt{1+y_0}\over2}\right)
- T\log(1+y_0).
$$
Using (5.11) a calculation gives
$$\eqalign{
r\log y_0 & = r\log {r\over T} + {r^2\over2T} +
O\left({r^4\over T^3}\right),\cr
- 2r\log\left({1+\sqrt{1+y_0}\over2}\right)& = -{r^2\over2T^2}
- {r^3\over16T^2} + O\left({r^4\over T^3}\right),\cr
- T\log(1+y_0) & = -r + {r^3\over24T^2} + O\left({r^4\over
T^3}\right),\cr}
$$
and the $O$-terms admit an asymptotic expansion in powers of $r/T$.
Therefore we obtain
$$
{\cal F}(y_0) - r\log4 = r\log\left({r\over4{\rm e}T}\right) + \sum_{j=3}^N
c_jr^jT^{1-j}  + O_N(r^{N+1}T^{-N})\leqno(5.12)
$$
for any given integer $N\ge3$ and some effectively computable
real constants $c_j\;(c_3 = -1/48)$. As
$$
{\cal F}''(y_0) \sim - {T^2\over r}\qquad(T\to\infty),\leqno(5.13)
$$
it follows that the dominant contribution to $I$ is a multiple of
$$
\leqalignno{\qquad&
T^{{1\over2}-2\tau}r^{2\tau-{3\over2}}
\exp\left\{-{\txt{1\over4}}G^2\log^2(1+y_0)
+i{\cal F}(y_0)-ir\log4\right\}&(5.14)\cr
=&T^{{1\over2}-2\tau}r^{2\tau-{3\over2}}\cr\times&
\exp\left\{-{\txt{1\over4}}G^2\log^2(1+y_0)
+ ir\log\left({r\over4{\rm e}T}\right) + i\sum_{j=3}^Nc_jr^jT^{1-j}
   + O_N(r^{N+1}T^{-N})\right\}.\cr}
$$
This is understood in the following sense: the remaining terms in
the evaluation of $I$ are either negligible, or similar in nature to
(5.14) (meaning that the oscillating exponential factor is the same,
which is crucial), only of the lower order of magnitude than (5.14).
We shall show now briefly show how the saddle point method does indeed
lead to this assertion.

\medskip
To see this we turn back to the integral in (5.9). We use the
techniques which were used in establishing (7.1.30) and (7.1.31) of
[Mo6]. With $y_0$ as in (5.11) we have that the relevant integral
is equal to
$$
y_0{\rm e}^{-\pi i/4}\int_{-\xi_0}^{\xi_0}f_0(\xi){\rm
e}^{if(\xi)}\d \xi\qquad (\xi_0 = r^{\eps-1/2}),\leqno(5.15)
$$
plus as error term which is $\ll_\eps \exp(-r^\eps)$. This error term
is negligible if
$$
r \,\ge\, (\log T)^{C(\eps)}\leqno(5.16)
$$
with $C(\eps)\;(>0)$ sufficiently large. The functions appearing
in (5.15) are ($\xi$ is the variable of integration)
$$
\leqalignno{
f_0(\xi) &:= y^{2\tau-3/2}\exp(-{\txt{1\over4}}G^2\log^2(1+y)),
\quad y := y_0 + y_0 \xi{\rm e}^{-\pi i/4},&(5.17)\cr
f(\xi) & := {\cal F}_-(\xi;r,T)
= r\log y - 2r\log\left({1+\sqrt{1+y}\over2}\right)
- T\log(1+y),\cr}
$$
where we assume that (5.1) and (5.16) hold. This enables us
to replace $f_0(\xi)$ with
$$
y_0^{2\tau-3/2}\exp(-{\txt{1\over4}}G^2\log^2(1+y_0)),
$$
on expanding $f_0(\xi)$ into its Taylor series at $y_0$. Likewise,
since ${\cal F}'(y_0) = 0$,
$$
f(\xi) = {\cal F}(y_0 +  y_0 \xi{\rm e}^{-\pi i/4})
= {\cal F}(y_0) + \hf iy_0^2\xi^2(-{\cal F}''(y_0)) + G(\xi;r,T,y_0),
$$
say, where $G$ can be expanded into Taylor series and
$$
G(\xi;r,T,y_0) \ll y_0^3\xi_0^3ry_0^{-3} = r^{3\eps-1/2}.
$$
After this the ensuing integrals are evaluated by using
the formula (proved by induction on $k$)
$$
\int_{-\Xi_0}^{\Xi_0}\xi^{2k}{\rm e}^{-{1\over2}c\xi^2}\d \xi
= 2^{k+{1\over2}}\G(k+\hf)c^{-{1\over2}-k} +
O_k(c^{-1}\Xi_0^{k-1}{\rm e}^{-{1\over2}c\Xi_0^2}),\leqno(5.18)
$$
provided that
$$
k = 0,1,2,\;\ldots\;,\;c>0,\;\Xi_0 > 0,\; \Xi_0\sqrt{c} \ge
1.\leqno(5.19)
$$
In our case
$$
c = - y_0^2F''(y_0) > 0,\; \Xi_0 = r^{\eps-1/2},\;
\Xi_0\sqrt{c} \asymp r^\eps,
$$
so that (5.19) is satisfied. Collecting all the estimates,
we see that the major contribution to $I$ is indeed furnished
by (5.14).

In the case when the integral in (5.9) has no saddle point,
i.e., the case of ${\cal F}_+$, we turn the segment of integration
by the angle $r^{-1/2}$, say, to obtain that the contribution
of the integral is in this case negligible.

In the case when (5.16) fails, more precisely when
$$
|r| \le (\log T)^{C(\eps)},
$$
we apply the technique of [Mo6, Lemma 5.2], to see that the integral
in question in the above range is $\ll T^{{1\over2}-2\tau}$, which
is sufficiently sharp for our purposes.

\bigskip
\centerline{\bb 6. The weighted fourth moment when
$\hf<\s<{3\over4}$}
\bigskip
With the use of Theorem 1 and the asymptotics of Section 5 we can
derive the explicit formula for the fundamental function
$$
\leqalignno{
I_2(T,\tau,G) := &{1\over\sqrt{\pi}G}\int_{-\infty}^\infty |\z(\tau
+ it + iT)|^4{\rm e}^{-(t/G)^2}\d t&(6.1)\cr
&\left(\hf<\tau<{\txt{3\over4}}, T^{1/3+\eps} \le G \le
T^{1-\eps}\right).
}
$$
This formula, as in the case when $\tau =
1/2$ (see [I2], [Mo6]), can be integrated over $T$. It
will then lead to  explicit results on
the function $E_2(T,\s)$, the error term in the asymptotic
formula for $\int_0^T|\z(\s+it)|^4\d t$.
Our result on $I_2(T,\s)$ and its integral is given by

\bigskip
THEOREM 3. {\it
If $I_2 (T, \sigma, G)$ is given by $(6.1)$,
$\hf < \s < {3\over4}, T^{1/3+\eps} \le G \le T^{1-\eps}$,
$Y_0 = (\k_j/T)(\sqrt{1+(\k_j/4T)^2} + \k_j/(2T)),$
then we have
$$
\leqalignno{
&I_2 (T, \sigma, G)  \sim O(1) \,+&(6.2)
\cr&
+ \;C(\sigma)T^{{1\over2} - 2\sigma}
\sum\limits_{\kappa_j
\leq TG^{-1} \log T} \alpha_j \kappa_j^{2 \sigma-3/2} H^2_j
(\hf) H_j (2 \sigma - \hf){\rm e}^{-{1\over4}G^2\log^2(1+Y_0)}\cr&
\times \sin \Bigl(\kappa_j \log
{\kappa_j\over 4{\rm e}T} +c_3\k_j^3T^{-2}\Bigr).\cr}
$$
We also have, for ${\bar Y}_0 =
(\k_j/ V)(\sqrt{1+(\k_j/4V)^2} + \k_j/(2V)),$ and
$V^{1/3+\eps}\leq
G \leq V^{1- \varepsilon}\;(D>0)$,
$$
\leqalignno{
&\int\limits^V_0 I_2 (T, \sigma ; G)\d T \sim
{\z^4(2\s)\over\z(4\s)}V +
{V\over3-4\s}{\left({V\over2\pi}\right)}^{2-4\s}{\z^4(2-2\s)\over\z(4
-4\s)}
&(6.3)\cr&
+ V^{2-2\s}(a_0(\s) + a_1(\s)\log V + a_2(\s)\log^2V) \cr&
+ \;C (\sigma) V^{{3\over2}-2 \sigma}
\sum\limits_{ \kappa_j \leq VG^{-1} \log T}
  \alpha_j \kappa_j^{2\sigma - 5 /2} H^2_j (\hf)H_j (2
\sigma - \hf){\rm e}^{-{1\over4}G^2\log^2(1+{\bar  Y}_0)}\cr&
\times\cos
\left(\kappa_j \log\bigl({\kappa_j\over4eV}\bigr) +c_3\k_j^3V^{-2}
\right) + O(G) + O(V^{1/3} \log^D V)
\cr}
$$
with suitable constants $C(\s), C_1(\s)$, and $a_j(\s)$,
which may be explicitly evaluated. The meaning of the
symbol $\sim$ is that besides the spectral sums in $(6.2)$-$(6.3)$ a finite
number of other sums are to appear, each of which is similar in nature
to the the corresponding sum above, but of a lower order
of magnitude.}

\bigskip
{\bf Proof of Theorem 3}.
The meaning of the symbol $\sim$ was already explained after
(5.14). Each of the omitted sums is either negligibly small, or
similar in structure to the ones appearing above, namely it has the same
oscillatory factors as the corresponding sums above. When estimated,
their contribution will be (by a power of $T$ or $V$) smaller than
the contribution  of the sums in (6.2) and (6.3).

To prove Theorem 3, we use (2.3)-(2.8) of Theorem 1.  The
derivation of (6.2) is similar to the proof of Theorem 5.2 of [I2]
or Theorem 5.1 of [Mo6], starting from the spectral
decomposition of ${\cal L}(g;\tau,\tau)$ when $\tau = 1/2$.
Thus we shall be relatively brief, noting that the sum in (6.2)
comes from the discrete spectral part (2.5) and (5.14). We shall need
(5.14)
with $\tau = \s, \,\hf < \s < 1$. The weight function $g$ will  be
(5.1), hence
$$
g_c(x) = {\rm e}^{-{1\over4}G^2x^2}\cos(xT).\leqno(6.4)
$$
In view of the expressions for $M_\ell^*\;(\ell = 0,\ldots,6)$
(see (3.74)--(3.89)) of the main term (cf. ${\cal Z}_r(\tau,g)$
in (2.4)) will be O(1), as will also be the contribution of
${\cal Z}_h(\tau,g)$ in (2.7). The contribution of ${\cal
Z}_c(\tau,g)$,
given by the integral in (2.6), is estimated by the use of
(5.14). It will be $O(1)$ plus the term which is
$$
\eqalign{&
\ll \log^2T\int_{-TG^{-1}\log T}^{TG^{-1}\log T}
|\zr|^4|\z(2\s-\hf)|^2T^{1/2-2\s}(|r|+1)^{2\s-3/2}\d r\cr&
\ll (TG^{-1})^{{5\over4}+2\s-{3\over2}}T^{{1\over2}-2\s}\log^CT\cr&
= T^{1/4}G^{1/4-2\s}\log^CT \le 1\cr}
$$
for $G \ge T^{1/3}$, since $\s > \hf$. Here
we used the trivial bound $1/|\zt| \ll \log|t|$, coupled with the
Cauchy-Schwarz inequality for integrals and the bounds (see [I1])
$$
\int_0^T|\zt|^8\d t \ll T^{3/2}\log^CT,\quad
\int_1^T|\z(\tau+it)|^4\d t \ll T \;\;(\hf < \tau \le 1).
$$

\medskip
To prove (6.3), we integrate first the spectral decomposition
of Theorem 1  from $V$ to $2V$, eventually replacing
$V$ by $V2^{-j}$  and summing over $j\in\NN$.
When we apply  (5.14) and integrate, we essentially have to integrate
$T^{{1\over2}-2\tau-ir}$ over $T$, which accounts for the
increase in order of $T/\k_j$ in (6.3), and one can check
that integration will transform the sine into cosine.
Here care should be exerted when one computes the main term
on the right-hand side of (6.3). This is given (cf.
$M({\rm p}_\tau;g)$) by eqs.(3.74)--(3.89).
In the evaluation we make repeatedly use of the formula
(see [I2, Section 5.1])
$$
{\G^{(k)}(s)\over\G(s)} =
\sum_{j=0}^kb_{j,k}(s)\log^js + c_{-1,k}s^{-1} +\ldots +c_{-r,k}s^{-r}
+ O_r(|s|^{-r-1})\leqno(6.5)
$$
for any fixed integers $k\ge1, r\ge0$, where each $b_{j,k} \,(\sim
b_{j,k}$
for a suitable constant $b_{j,k}$) has an asymptotic expansion in
non-positive powers of $s$. It transpires that one encounters integrals
of the type
$$
\leqalignno{&
{1\over\sqrt{\pi G}}
\int_{-\infty}^\infty \log^r(\hf+iT+it){\rm e}^{-(t/G)^2}\d t&(6.6)\cr&
= {1\over\sqrt{\pi}}
\int_{-\infty}^\infty \log^r(\hf+iT+iuG){\rm e}^{-u^2}\d u\cr&
= {1\over\sqrt{\pi}}
\int_{-\log T}^{\log T} \log^r(\hf+iT+iuG){\rm e}^{-u^2}\d u
+ O_A(T^{-A}),\cr}
$$
for any fixed $A>0$. For $|u| \le \log T$ one has
the power series expansion
$$\eqalign{&
\log^r(\hf+iT+iuG) = \log^r(iT)\cr&
+ \sum_{k=1}^r{r\choose k}(\log iT)^{r-k}\left(
{uG\over T} + {1\over2iT} - {1\over2}\left({uG\over T} +
{1\over2iT}\right)^2
+ \ldots\right)^k\cr}
$$
which is inserted in (6.6). The evaluation is completed by applying
(5.18). The main term in (6.3) is the same one as in (6.2) (with
$V$ replacing $T$), and the constant standing in front of the term
$T^{3-4\s}$ was first explicitly evaluated by Ka\v c\.enas [K1, K2].
The contributions of ${\cal Z}_h(\tau,g) $ and ${\cal Z}_c(\tau,g)$
will be absorbed by the error terms after integration. This ends our
discussion
of Theorem 3.

\medskip
Next, we consider $E_2(T,{3\over4})$ by using (3.88)--(3.90). We
obtain,
with suitable constants $A_j$, which may be explicitly evaluated,
$$
\leqalignno{\int_0^T|\z({\txt{3\over4}}+it)|^4\d t &=
{\z^4({\txt{3\over2}})\over\z(3)}T &(6.7)\cr&+ T^{1/2}(A_0  + A_1\log T
+
A_2\log^2T) + E_2(T,{\txt{3\over4}})\cr}
$$
with
$$
E_2(T,{\txt{3\over4}}) \,\ll\, T^{1/2}\log^3T.\leqno(6.8)
$$
Note that the bound (6.8) for the error term is, by a log-factor,
larger than the order of the second main term in (6.7). Indeed, it
is very plausible that the bound (6.8) is far from the truth and
that we have
$$
E_2(T,\s) \,\ll_\eps\, T^{3/2-2\s+\eps}      \qquad(\hf < \s <
{\txt{3\over4}})\leqno(6.9)
$$
and
$$
E_2(T,\s) \,\ll_\eps\, T^{\eps}     \qquad(
{\txt{3\over4}} \le \s < 1).\leqno(6.10)
$$
Here and later $\eps$ denotes arbitrarily small, positive constants,
not necessarily the same ones at each occurrence, and $f \ll_\eps g$
means that the $\ll$--constant depends on $\eps$. Also note that $C$
will denote a generic positive constant.

The conjectures (6.9)-(6.10) were made in [I6]. They are very
strong, since they imply that $\zt \ll_\eps |t|^{1/8+\eps}$
and $\z(\s+it) \ll_\eps |t|^\eps$ for $\s \ge {3\over4}$.
They are the analogues of the conjectures for the true order of the
error term $E_1(T,\s)$ in (1.5) (see [Ma]).
What seems possible to prove at present
  for the range    ${\txt{3\over4}} \le \s < 1$
is (cf. [I6, Th. 2])
$$
\int_0^T|\z({\txt{3\over4}}+it)|^4\d t =
{\z^4(2\s)\over\z(4\s)}T + O(T^{2-2\s}\log^3T),\leqno(6.11)
$$
which is far from the conjectured bound (6.10). In view of (6.11)
there seems to be no point in further estimation of $E_2(T,\s)$
when ${3\over4} < \s < 1$, since the bounds that seem
obtainable from the spectral decomposition are weaker than (6.11).
When $\s = 1$ we have (see [I3])
$$
\int_1^T|\z(1+it)|
^4\d t = {\z^4(2)\over\z(4)}T + O(\log^4T),
$$
so that this case is covered, too (${\z^4(2)/\z(4)} = \pi^2/72$).

%\break

\bigskip
\centerline{\bb 7. The fourth moment when $\hf < \s < {\txt{3\over4}}$}
\bigskip
We have prepared the groundwork for the results on $E_2(T,\s)$,
the error term for the fourth moment off the critical line (see
(1.6)), in the previous sections. Now we can proceed with the
statement of our results.

\bigskip
{THEOREM 4}. {\it If $\s$ is a fixed number such that $\hf < \s <
{3\over4}$, and $E_2(T,\s)$ is defined by} (1.6),
{\it then with suitable constants $a_j(\s)$ we have}
$$
\leqalignno{\int_0^T|\z(\s+it)|^4\d t &=
{\z^4(2\s)\over\z(4\s)}T +
{T\over3-4\s}{\left({T\over2\pi}\right)}^{2-4\s}{\z^4(2-2\s)\over\z(4
-4\s)}
&(7.1)\cr&
+ T^{2-2\s}(a_0(\s) + a_1(\s)\log T + a_2(\s)\log^2T) +
E_2(T,\s),\cr}
$$
{\it where with some } $C>0$
$$
E_2(T,\s) \,\ll\, T^{2/(1+4\s)}\log^CT      \qquad(\hf < \s <
{\txt{3\over4}}).\leqno(7.2)
$$
{\it Moreover,}
$$
E_2(T,\s) \,=\, \Omega_\pm(T^{{3\over2}-2\s})      \qquad(\hf < \s <
{\txt{3\over4}}).\leqno(7.3)
$$
{\it More precisely, there exist constants $A = A(\s) > 1$ and
$B = B(\s) > 0$ such that, for $T \ge T_0(\s)$, every interval
$[T,\,AT]$ contains points $t_1 = t_1(\s)$ and $t_2 = t_2(\s)$ such
that}
$$
E_2(t_1,\s) > Bt_1^{{3\over2}-2\s},\quad E_2(t_1,\s) <
-Bt_1^{{3\over2}-2\s}
  \qquad(\hf < \s < {\txt{3\over4}}).\leqno(7.4)
$$
\bigskip
{\bf Remarks}. As usual, $f(x) = \Omega_\pm(g(x))$ means that
we have $\limsup_{x\to \infty}
f(x)/g(x) > 0$ and $\liminf_{x\to\infty} f(x)/g(x) < 0$ for a given
$g(x) > 0 \,(x\ge x_0)$.
Note that $3-4\s > 2/(1+4\s)$ for $\hf < \s < {1+\sqrt{2}\over4}$
and that $2-2\s > 2/(1+4\s)$ for $\s < {3\over4}$.
  Thus our bound for the error term $E_2(T,\s)$
is already larger than the second main term in (1.6) unless
$\hf < \s < {1+\sqrt{2}\over4}$, but the bound in question is probably
much too large (recall the conjectural bounds
(6.9)--(6.10) for the order of $E_2(T,\s)$).

\bigskip
THEOREM 5. {\it Let $E_2(T,\s)$ be given by} (1.6).
{\it If $\s$ is a fixed number such that $\hf < \s <
{3\over4}$, then for suitable $C = C(\s)>0$  we have
$$
\int_0^T|E_2(t,\s)|^{4\s}\d t \ll T^{2}\log^CT.\leqno(7.5)
$$
We also have, for any constant $A\ge 1$,}
$$
\int_0^T|E_2(t,\s)|^{A}\d t  \gg T^{1+A({3\over2}-2\s)}.\leqno(7.6)
$$

\medskip
Note that when $\s = \hf$, (7.5) reduces to
$$
\int_0^TE^2_2(t)\d t \ll T^{2}\log^CT,\leqno(7.7)
$$
where $E_2(T) = E_2(T,\hf)$ is the error term in the formula for
the fourth moment of $|\zt|$. The bound (7.7)
is the sharpest one known (see [IM2], [Mo6]) and essentially best
possible,
since we have (see [I7])
$$
\int_0^TE^2_2(t)\d t \gg T^2.\leqno(7.8)
$$
The lower bound in (7.6), when $A=2,\,\s = \hf$, reduces
to (7.7).
Note  that the conjecture (6.9) would furnish the upper bound
$$
\int_0^T|E_2(t,\s)|^{A}\d t  \ll_\eps T^{1+A({3\over2}-2\s)+\eps},
\leqno(7.9)
$$
which is (up to `$\eps$') the same as the lower bound (7.6). The upper
bound in (7.5), on the other hand, is much weaker than (6.9). This
reflects, in general, the situation with $E_2(T,\s)$: as $\s$ increases
from $1\over2$ to $3\over4$, the quality of the bounds (either pointwise
or in the mean square sense) decreases. The same phenomenon also occurs
with bounds for the mean square of $\z(s)$ off the critical line (see
[Ma]).

Finally we remark that, by H\"older's inequality for integrals,
(7.5) implies the mean square bound
$$
\int_0^TE^2_2(t,\s)\d t \ll T^{1+ 1/(2\s)}\log^CT,\leqno(7.10)
$$
which may be compared to (7.7).

\medskip
\centerline{\bb 8. Proof of the bounds when $\hf < \s < {3\over4}$}
\medskip
In this section we shall prove Theorem 4 and Theorem 5,
formulated  in the preceding section. First we prove
(7.1)-(7.2). We rewrite (6.3) of Theorem 3 as
$$
\int_0^T I_2(t,\s;G)\d t = M(T,\s) + S(T,\s;G) + R(T,\s;G),\leqno(8.1)
$$
say, where the main term is
$$
\leqalignno{M(T,\s) &:=
{\z^4(2\s)\over\z(4\s)}T +
{T\over3-4\s}{\left({T\over2\pi}\right)}^{2-4\s}{\z^4(2-2\s)\over\z(4
-4\s)}
&(8.2)\cr&
+ T^{2-2\s}(a_0(\s) + a_1(\s)\log T + a_2(\s)\log^2T),\cr
}
$$
$$
\leqalignno{
S(T,\s;G) &:=
C (\sigma) T^{{3\over2}-2 \sigma}
\sum\limits_{j=1}^\infty
  \alpha_j \kappa_j^{2\sigma - 5 /2} H^2_j (\hf)H_j (2
\sigma - \hf){\rm e}^{-{1\over4}G^2\log^2(1+Y_0)}&(8.3)\cr& \times\cos
\left(\kappa_j \log\bigl({\kappa_j\over4eT}\bigr) +c_3\k_j^3T^{-2}
\right)\cr}
$$
is the spectral part, and the rest (error term) is
$$
R(T,\s;G) := O(G) + O(T^{1/3}\log^DT).\leqno(8.4)
$$
We suppose that $T^{1/3+\eps} \le G = G(T) \le T^{1-\eps}$ and put
first in (8.1)
$$
T_1 = T - G\log T,\; T_2 = 2T + G\log T.
$$
Then
$$
\eqalign{&
\int_{T_1}^{T_2}I_2(t,\s;G)\d t = \int_{-\infty}^{\infty}
|\z(\s+iu)|^4\left({1\over\sqrt{\pi}G}\int_{T_1}^{T_2}
{\rm e}^{-(t-u)^2/G^2}\d t\right)\d u\cr&
\ge \int_{T}^{2T}\z(\s+iu)|^4\left({1\over\sqrt{\pi}G}
\int_{T-G\log T}^{2T+G\log T}{\rm e}^{-(t-u)^2/G^2}\d t\right)\d u.\cr}
$$
But for $T\le u\le2T$ we have, by the change of variable $t-u = Gv$,
$$
\eqalign{&
{1\over\sqrt{\pi}G}\int_{T-G\log T}^{2T+G\log T}{\rm
e}^{-(t-u)^2/G^2}\d t
= {1\over\sqrt{\pi}}
\int_{(T-u)/G-\log T}^{(2T+u)/G\log T}{\rm e}^{-v^2}\d v\cr&
= {1\over\sqrt{\pi}}\int_{-\infty}^{\infty}{\rm e}^{-v^2}\d v
+ O\left(\int_{\log T}^\infty {\rm e}^{-v^2}\d v\right)
+ O\left(\int^{-\log T}_{-\infty} {\rm e}^{-v^2}\d v\right)\cr&
= 1 + O({\rm e}^{-\log^2T}),\cr}
$$
since $t - u\le0,\;2T - u\ge0$. Therefore, by (8.1) and the mean
value theorem, we obtain
$$
\leqalignno{&
\int_{T}^{2T}|\z(\s+it)|^4\d t \le
\int_{T_1}^{T_2}I_2(t,\s;G)\d t  + O(1)&(8.5)\cr&
= M(2T,\s) - M(T,\s) + O(G)\cr&
+ S(2T + G\log T,\s;G) - S(T - G\log T,\s;G)\cr&
+ R(2T + G\log T,\s;G) - R(T - G\log T,\s;G).\cr}
$$
A lower bound of a similar type for the first integral in (8.5) follows
by the same procedure if we take
$$
T_1 = T + G\log T,\; T_2 = 2T - G\log T.
$$
Putting together the bounds we obtain the following lemma, which is the
analogue of [I2, Lemma 5.1] or [IM3, Lemma 3].

\medskip
{\bf Lemma 11.} {\it With the notation introduced in} (8.1)-(8.2) {\it
and, for $T^{1/3+\eps} \le G \le T^{1-\eps}$, we have}
$$
\leqalignno{&
E_2(2T,\s) - E_2(T,\s)&(8.6)\cr&
\ll |S(2T + G\log T,\s;G)| + |S(2T - G\log T,\s;G)|\cr&
+ |S(T + G\log T,\s;G)| + |S(T - G\log T,\s;G)|\cr&
+ O(G) + O(T^{1/3}\log^DT).\cr}
$$

To return to the proof of (7.1)-(7.2), note that the $S$-sums can
be truncated at\break $TG^{-1}\log T$ with a
negligible error.  We estimate the exponential factors trivially, and
then use the
Cauchy-Schwarz inequality and the bound (4.2). Thus we have
$$
\eqalign{&
|S(2T + G\log T,\s;G)| + |S(2T - G\log T,\s;G)|\cr&
+ |S(T + G\log T,\s;G)| + |S(T - G\log T,\s;G)|\cr&
\ll T^{{3\over2}-2\s}(TG^{-1}\log T)^{2\s-{1\over2}}\log^CT\cr&
\le TG^{{1\over2}-2\s}\log^{C+1}T.\cr}
$$
This gives, by Lemma 11,
$$
E_2(2T,\s) - E_2(T,\s) \ll (TG^{{1\over2}-2\s} + G + T^{1/3})\log^CT
\ll T^{2\over1+4\s}\log^CT\leqno(8.7)
$$
with the choice
$$
G \=\ T^{2\over1+4\s}.
$$
 From (8.7) the bound (7.2) easily follows. Note that an explicit
value $C = C(\s)$ in Theorem 4 can also be worked out without trouble.

\medskip
To prove the omega result (7.3) we argue similarly as in the case
of the proof of the omega-result (see [I7], [Mo4], [Mo6])
$$
E_2(T) \;=\;\Omega_\pm(T^{1/2}).\leqno(8.8)
$$
Instead of the (modified) Mellin transform
$$
{\cal Z}_2(s) := \int_1^\infty |\zx|^4x^{-s}\d x\qquad(\Re s > 1)
\leqno(8.9)
$$
used for the proof of (8.8), we need to use the function
$$
{\cal Z}_{2}(s,\tau) := \int_1^\infty
|\z(\tau+ix)|^4x^{-s}\d x\qquad(\hf < \tau < 1,\;\Re s > 1).
$$
The spectral decomposition of ${\cal Z}_{2}(s,\tau) $ is effected much
in the same way as was the spectral decomposition of ${\cal Z}_2(s)$
(see [Mo4], [Mo6]). The major difference relevant for the omega results
is that, in the case of ${\cal Z}_2(s)$ the simple poles are located
at $s=\hf \pm i\k_j$, while in the case of ${\cal Z}_{2}(s,\tau) $
the simple poles are located at $s= 3/2-2\tau \pm i\k_j$. Hence,
instead of (8.8), we obtain the omega result (7.3).
In the course of the proof one needs the non-vanishing of
$L_\s(\k)$ for infinitely many $\k$ (see (4.1)), which follows
from (4.3) of Theorem 2. The function
${\cal Z}_{2}(s,\tau) $ admits meromorphic continuation over $\CC$
where,
unless $s$ lies in a neighborhood of its pole, it is of polynomial
growth
in $|s|$ for a fixed $\s$. This follows analogously as in [Mo6] and
[IJM]. The crucial result is analogue of Lemma 2 of [I7], which
in this case will imply that
$$
\leqalignno{&
\int_0^\infty \int_0^tE_2(u,\s)\d u\cdot{\rm e}^{-t/T}\d t
&(8.10)\cr&
\sim T^{{7\over2}-2\s}\Re\left\{\sum_{j=1}^\infty \a_jH_j^2(\hf)
H_j(2\s-\hf)R_{1,\s}(\k_j)\right\}\quad(T\to\infty),\cr}
$$
where $R_{1,\s}(\k_j) \ll_\eps \exp(-(\hf\pi-\eps)\k_j)$.
 From (8.10) we obtain (7.4) with the aid of [I7, Lemma 3].
With (7.4) at our disposal, we prove easily (7.3).
Let $t_1$ be as in (7.4). Then
$$
At_1^{{5\over2}-2\s} < \int_0^{t_1}E_2(t,\s)\d t
\le \left(\int_0^{t_1}|E_2(t,\s)|^a\d t\right)^{1/a}t_1^{(a-1)/a}
$$
for $a > 1$ by H\"older's inequality, and for $a=1$ it is trivial. In
view
of $T \le t_1 \le BT$ we obtain
$$
A^aT^{1+a({3\over2}-2\s)} \le A^at_1^{1+a({3\over2}-2\s)}
\le \int_0^{t_1}|E_2(t,\s)|^a\d t
\le \int_0^{BT}|E_2(t,\s)|^a\d t.
$$
Changing $T$ to $T/B$ we obtain (7.3).

\medskip

It remains to prove (7.5) of Theorem 5. We shall follow the proof of
[Mo6, Theorem 5.3], making the necessary modifications. We wish to
obtain an upper bound for $R$,  the number
of well-spaced points $\{t_r\}\;(r = 1,\ldots,R)$ for which
$E_2(t_r,\s) \ge V > 0$ (the case when $E_2(t_r,\s) \le -V$ is
analogously treated, so we may consider only the former case),
where
$$
\leqalignno{
T \le t_1 < \cdots &< t_R\le 2T,\quad
t_{r+1}-t_r \ge V\log^{-C-1}T,&(8.11)\cr
& T^{1\over4\s}\log^{C_2}T \le V \le
T^{2\over1+4\s}\log^{C_3}T,
}
$$
for suitable $C_j > 0$.
We put
$$
U = 2^{-\ell}t_r \quad(\ell = 1,\ldots,L),\quad G =
V\log^{-C_4}T,\quad 2^{-L}T \asymp T^{(4\s+1)/(8\s)},
$$
which gives
$$
E_2(t_r) = \sum_{\ell=1}^L\sum_{r=1}^R\left\{E_2(2^{1-\ell}t_r,\s)
- E_2(2^{-\ell}t_r,\s)\right\} + O(T^{1/(4\s)}\log^CT)
$$
by (7.2). Therefore we obtain
$$ \hf RV \le
\sum_{\ell=1}^L\sum_{r=1}^R\left\{E_2(2^{1-\ell}t_r,\s) -
E_2(2^{-\ell}t_r,\s)\right\},\leqno(8.12)
$$
and we now apply Lemma 11. We may truncate each sum over
$\k_j$ so that $\k_j \le TG^{-1}\log T$, and also expand into
Taylor series the factor
$$ \exp\left(ic_3\k_j^3(U\pm G\log T)^{-2}\right)
$$
and higher power exponentials coming from (5.14), noting that
the main contribution will come from the constant term, namely
unity. This is important, since this procedure allows us to
relax the condition $G \ge V^{1/3+\eps}$ in (6.3) in such a way
that $G$ and $U$ lie in a permissible range.
Instead of $W(K,\ell;z)$ (cf. [Mo6, eq.\ (7.2.19)]) we have
($\tau(r,\ell) = 2^{1-\ell} + G\log T,\,\Re z = 1/\log T$) now
$$
\eqalign{ W(K,\ell;z) &:= \sum_{K<\k_j\le
2K}\a_jH_j^2(\hf)|H_j(2\s-\hf)|\k_j^{2\s-1}
\left|\sum_{r=1}^R\tau(r,\ell)^{{3\over2}-2\s+z+i\k_j}\right|\cr&
\ll K^{2\s-1}\sum_{K<\k_j\le 2K}\a_jH_j^2(\hf)|H_j(2\s-\hf)|
\left|\sum_{r=1}^R\tau(r,\ell)^{{3\over2}-2\s+z+i\k_j}\right|.\cr}
$$
To the
sum over $\k_j$ we apply the Cauchy-Schwarz inequality, noting
that for $\s > \hf$
$$
\eqalign{& \sum_{K<\k_j\le
2K}\a_j|H_j(\hf)H_j(2\s-\hf)|^2\cr& \le \left(\sum_{K<\k_j\le
2K}\a_jH_j^4(\hf) \sum_{K<\k_j\le
2K}\a_jH_j^4(2\s-\hf)\right)^{1/2}\cr& \ll K^2\log^CK,\cr}
$$
since both sums above are bounded by $K^2\log^CK$. For the sum
with $H_j^4(\hf)$ this is [Mo6, Theorem 3.4], and the other sum
is treated analogously. This yields
$$
W^2(K,\ell;z) \ll
K^{4\s}\log^CK\sum_{K<\k_j\le 2K}\a_jH_j^2(\hf)
\Bigl|\sum_{r=1}^R\tau(r,\ell)^{{3\over2}
-2\s+z+i\k_j}\Bigr|^2.\leqno(8.13)
$$
With (8.13) we obtain, on applying [Mo6, eq.\ (5.6.3)]
(this is a variant of the  spectral large sieve), the uniform bound
$$
W^2(K,\ell;z) \ll K^{4\s+1}\log^CK(K+TV^{-1})RT^{3-4\s}2^{-\ell}.
$$
With the aid of (8.12) this yields, similarly as in [Mo6, Chapter 5]
$$
\eqalign{ RV& \ll
\max_{K\le TG^{-1}\log T} K^{-3/2}(K^{2\s+1} +
T^{1/2}V^{-1/2}K^{2\s+{1\over2}})R^{1/2}T^{3/2-2\s}\log^CT\cr& \ll
R^{1/2}TV^{{1\over2}-2\s}\log^CT.\cr}
$$
Therefore we obtain
$$
R \ll T^2V^{-1-4\s}\log^CT,\leqno(8.14)
$$
which easily leads to
(7.2). The part where $|E_2(t,\s)| \le T^{1/(4\s)}\log^CT$ is
trivial, so we may restrict integration to the set $\cal S$, where
$|E_2(t,\s)| \ge T^{1/(4\s)}\log^CT$. Consider the subset ${\cal
S}_V$ of $\cal S$, where $V \le |E_2(t,\s)| < 2V, t\in {\cal S}
\cap [\hf T, T]$. We divide the interval $[\hf T,\,T]$ into
subintervals of length $V\log^CT$, allowing the end subintervals
to be possibly shorter. Then the number $R = R_V$ of those subintervals
(considering separately subintervals with even and odd indices)
which contain a point from ${\cal S}_V$ is bounded by (8.14).
Hence we have
$$
\int_{{\cal S}_V} |E_2(t,\s)|^{4\s} \d t
\ll R_VV\log^CTV^{4\s} \ll
T^{2}\log^CT,
$$
and since there are $O(\log T)$ choices for $V$, we have
$$
\int_{\hf T}^T|E_2(t,\s)|^{4\s}\d t \ll T^{2}\log^CT.
$$
Replacing $T$ by $T2^{-j}$ and summing the above bounds over $j\in\NN$
we obtain (7.5). The proof of Theorem 5 is complete.

\vfill
\eject

\topglue 2cm
\centerline{{\bf REFERENCES}}
\vskip1cm
\item{[BM]} R.W. Bruggeman and Y. Motohashi, `A new approach to the
spectral theory of the fourth moment of the Riemann zeta-function',
To appear in  {\it Journal reine angew.\ Math.}

\item {[I1]} { A. Ivi\'c},  `The Riemann zeta-function', {\it John
Wiley
and Sons}, New York, 1985 (2n ed., Dover, 2003).

\item {[I2]} { A. Ivi\'c},  `Mean values of the Riemann zeta-function',
LN's {\bf82}, {\it Tata Institute of Fundamental Research}, Bombay,
1991
(distr. by Springer Verlag, Berlin etc.).

\item {[I3]} A. Ivi\'c, `The moments of the zeta-function on the
line $\s = 1$', {\it Nagoya Math.\ J}. {\bf135}(1994), 113-120.

\item{ [I4]} { A. Ivi\'c},  `On the fourth moment of the Riemann
zeta-function', {\it Publs.\ Inst.\ Math. (Belgrade)}
{\bf57(71)}(1995), 101-110.

\item{[I5]} { A. Ivi\'c},  `The Mellin transform and
the Riemann zeta-function', in
{\it ``Proceedings of the Conference on Elementary and
Analytic Number Theory
(Vienna, July 18-20, 1996)"}, Universit\"at Wien \& Universit\"at f\"ur
Bodenkultur, Eds. W.G. Nowak and J. Schoi{\ss}enge-\break ier,
Vienna 1996, 112-127.

\item {[I6]} A. Ivi\'c,  Some problems on mean values of the Riemann
zeta-function, {\it Journal de Th\'eorie des Nombres Bordeaux} {\bf8} (1996),
101-122.

\item{[I7]} A. Ivi\'c, `On the error term for the fourth moment of the
Riemann zeta-function', {\it J. London Math.\ Soc.}
{\bf60}(2)(1999), 21-32.

\item{[I8]} A. Ivi\'c,  `On the integral of the error term in
the fourth moment of the Riemann zeta-function', {\it Functiones
et Approximatio} {\bf28}(2000), 37-48.

\item{[I9]} A. Ivi\'c, `On the moments of Hecke series at central
points',
{\it Functiones et Approximatio} {\bf30}(2002), 49-82.

\item{ [IM1]}  A. Ivi\'c and Y. Motohashi,  `A note on the mean value
of
the zeta and L-functions VII', {\it  Proc.\ Japan Acad.\ Ser. A}
{\bf  66}(1990), 150-152.

\item{ [IM2]}{ A. Ivi\'c and Y. Motohashi}, `The mean square of the
error term for the fourth moment of the zeta-function',
{\it Proc.\ London Math. Soc.} (3){\bf66}(1994), 309-329.

\item {[IM3]} { A. Ivi\'c and Y. Motohashi},  `The fourth moment of the
Riemann zeta-function', {\it J. Number Theory} {\bf 51}(1995), 16-45.

\item {[IJM]}  A. Ivi\'c, M. Jutila and Y. Motohashi,
`The Mellin transform of powers of the zeta-function', {\it Acta
Arith.} {\bf95}(2000), 305-342.

\item{[K1]} A. Ka\v c\.enas, `Mean values of the Riemann zeta-function
in the critical strip', {\it Doctoral Thesis},
Vilnius University, Vilnius, 1996.

\item{[K2]} A. Ka\v c\.enas, `Asymptotics of the fourth power moment
of the Riemann zeta-function in the critical strip',
{\it Lithuanian Math. J.} {\bf36}(1996), 32-44.

\item {[KS]} S. Katok and P. Sarnak, `Heegner points, cycles and Maass
forms', {\it Israel J. Math.} {\bf84}(1993), 193-227.

\item{[L]} N.N. Lebedev, `Special functions and their applications',
{\it Dover Publications, Inc.}, New York, 1972.

\item {[Ma]} K. Matsumoto, `Recent developments in the mean square
theory of the Riemann zeta and allied functions', in "{\it Number
Theory"}, eds. R.P. Bambah et al., Hindustan Book Agency \&
Indian National Science Academy, Birkh\"auser, 2000, pp. 241-286.

\item {[Mo1]} Y. Motohashi, `Spectral mean values of Maass wave
form $L$-functions', {\it J. Number Theory} {\bf42}(1992),
258-284.

\item {[Mo2]} Y. Motohashi, `An explicit formula for the fourth
power mean of the Riemann zeta-function', {\it Acta
Math.\/} {\bf170}(1993), 181-220.

\item {[Mo3]} Y. Motohashi, `The mean square of Hecke $L$-functions
attached to holomorphic cusp forms', {\it Kokyuroku Res.\ Inst.\ Math.
Kyoto Univ.} {\bf886}(1994), 214-227.

\item {[Mo4]} Y. Motohashi,  `A relation  between the Riemann
zeta-function and the hyperbolic Laplacian', {\it Annali Scuola Norm.\
Sup.\ Pisa, Cl.\ Sci.\ IV ser.} {\bf 22}(1995), 299-313.

\item {[Mo5]} Y. Motohashi,  `The Riemann zeta-function and the
non-Euclidean Laplacian', {\it Sugaku Expositions}, AMS {\bf 8}(1995),
59-87.

\item {[Mo6]} Y. Motohashi,  `Spectral theory of the Riemann
zeta-function', {\it Cambridge University Press}, Cambridge, 1997.

\item {[Mo7]} Y. Motohashi, `A note on the mean value of the zeta
and $L$-functions.\ XIV', {\it Proc.\ Japan Acad.} {\bf80A} (2004),
28-33.

\item{[WW]} E.T. Whittaker and G.N. Watson, `A Course of Modern
Analysis', {\it Cambridge University Press} (4th. ed.), London,
1963.

\bigskip
\leftline{\sevenrm Aleksandar Ivi\'c}
\leftline{\sevenrm Katedra Matematike RGF-a}
\leftline{\sevenrm Universiteta u Beogradu}
\leftline{\sevenrm Dju\v sina 7, 11000 Beograd}
\leftline{\sevenrm Serbia and Montenegro, ivic@rgf.bg.ac.yu}

\bigskip
\leftline{\sevenrm Yoichi Motohashi}
\leftline{\sevenrm Department of Mathematics}
\leftline{\sevenrm College of Science and Technology}
\leftline{\sevenrm Nihon University, Surugadai, Tokyo 101-8308}
\leftline{\sevenrm Japan, ymoto@math.cst.nihon-u.ac.jp;
http://www.ne.jp/asahi/zeta/motohashi/}

\bye